\documentclass{article}

\usepackage[doublespacing]{setspace}
\usepackage[fontsize=12pt]{fontsize}
\usepackage{amsthm, amsmath, amssymb}
\usepackage{comment}
\usepackage{setspace}
\usepackage{multirow}
\usepackage{enumitem}
\usepackage{xcolor}
\usepackage[utf8]{inputenc}
\usepackage{natbib}
\usepackage{url} % not crucial - just used below for the URL 
\usepackage{comment} % for commenting out section
\usepackage{natbib}
\usepackage{bbm}
\usepackage{bbold}
\usepackage{pifont}
\usepackage{amssymb}
\usepackage{mathbbol}
\usepackage{float}
\usepackage{algpseudocode}
\usepackage{algorithm}
\usepackage{float}
\usepackage{xcolor}
\usepackage{tabularray}
\usepackage{bm} % Math packages
\usepackage{setspace}
\usepackage{multirow}
\usepackage{lscape}
\usepackage{rotating}
\usepackage{adjustbox}
\usepackage{arydshln}
\usepackage{booktabs}
\usepackage{ragged2e}
\usepackage{tikz}
\usepackage[affil-it]{authblk}
\usepackage{xparse}
\usepackage{xr}
\usepackage{geometry} 
\usepackage[caption = false]{subfig} 
\usepackage{float}
\def\*#1{\mathbf{#1}}
\def\^#1{\mathbf{#1}}
\def\##1{\mathbb{#1}}

\DeclareSymbolFontAlphabet{\amsmathbb}{AMSb}%

\newtheorem{assumption}{Assumption}
\newtheorem{definition}{Definition}
\newtheorem{theorem}{Theorem}

\newtheorem{corollary}{Corollary}

\definecolor{blue}{RGB}{0,0,255}
\definecolor{slightgreen}{RGB}{34,139,34}
\definecolor{skyblue}{RGB}{30,180,255}
\definecolor{plum}{rgb}{0.56, 0.27, 0.52}

\def\A{\bm{A}}

\def\E{\mathbbm{E}}

\def\I{\bm{\mathcal{I}}} %\def\rI{{\bm{\cal{I}}}}
 
\def\J{\bm{\mathcal{J}}}

\def\O{\mathcal{O}}
   
\def\Q{\mathcal{Q}}

\def\T{\bm{T}} 
 
\def\bmV{\bm{V}} 
\def\V{\mathcal{V}}

\def\X{\bm{X}} \def\x{\bm{x}} 
\def\Y{\bm{Y}} \def\y{\bm{y}}  
\def\Z{\bm{Z}}

\newcommand\lowero{
  \mathchoice
    {{\scriptstyle\mathcal{O}}}% \displaystyle
    {{\scriptstyle\mathcal{O}}}% \textstyle
    {{\scriptscriptstyle\mathcal{O}}}% \scriptstyle
    {\scalebox{.7}{$\scriptscriptstyle\mathcal{O}$}}%\scriptscriptstyle
  }

\def\1{\bm{1}}  

\def\balpha{\bm{\alpha}}

\def\bbeta{\bm{\beta}}

\def\btheta{\bm{\theta}}  		 
\def\bmu{\bm{\mu}}

\def\bSigma{\bm{\Sigma}}

\def\bpi{\bm{\pi}}
\def\bzeta{\bm{\zeta}}

\def\bTheta{\boldsymbol{\Theta}}
\def\pconverge{\overset{p}{\rightarrow}}
\def\dconverge{\overset{d}{\rightarrow}}

\newcommand{\argmax}{\mathop{\mathrm{argmax}}} 
\newcommand{\argmin}{\mathop{\mathrm{argmin}}}

\title{Maximum Ideal Likelihood Estimation: A Unified Inference Framework for Latent Variable Models}

\author[1]{Jake Yizhou Cai}
\author[1]{Ting Fung Ma}

\affil[1]{Department of Statistics, University of South Carolina\\ 1523 Greene Street, SC 29208, United States}
\date{\today}

\begin{document}
\maketitle

\begin{abstract}
    This paper develops a unified estimation framework, the Maximum Ideal Likelihood Estimation (MILE), for general parametric models with latent variables. Unlike traditional approaches relying on the marginal likelihood of the observed data, MILE directly exploits the joint distribution of the complete data by treating the latent variables as parameters (the ideal likelihood). Borrowing strength from optimisation techniques and algorithms, MILE is a broadly applicable framework in case that traditional methods fail, such as when the marginal likelihood has non-finite expectations. MILE offers a flexible and robust alternative to established techniques, including the Expectation-Maximisation algorithm and Markov chain Monte Carlo. We facilitate statistical inference of MILE on consistency, asymptotic distribution, and equivalence to the Maximum Likelihood Estimation, under some mild conditions. Extensive simulations illustrative real-data applications illustrate the empirical advantages of MILE, outperforming existing methods on computational feasibility and scalability.
\end{abstract}


% \boxedtext{
% \begin{enumerate}[label=(\alph*)~]
% \item Key boxed text here.
% \item Key boxed text here.
% \item Key boxed text here.
% \end{enumerate}}

\maketitle

\section{Introduction}
\label{Introduction}
With the growth of data complexity, dependence are increasingly embedded into statistical models. Latent variable models, because of convenience to understand and interpret, are widely applied to dependent scenarios in domain ranging from causal inference \citep{causal2009} and spatio-temporal models \citep{INLA2009, Cressie2011} to advances in large language models \citep{LLM2023}. Traditionally, the Expectation-Maximisation (EM) algorithm \citep{EMHartley1958, EM1977} and its modifications (e.g., \cite{stEM2000, MCEM2001, MCEMreview2024}) remains the widely-used techniques to address information loss and latent structure. Collectively called EM-type algorithms, these methods provide perspectives across diverse settings, including mixture models and missing data problems. Maximising the conditional expectation of log-likelihood, EM-type algorithms yield to the maximum likelihood estimator (MLE), which serve as basis of inference. See \cite{EMbook2008} for reviews. Alternatively, Markov chain Monte Carlo (MCMC) \citep{MCMCbook2011, BDA3} offers a general framework under Bayesian settings, approximating the posterior distribution through sampling. Owing to the stability in case of complicated and implicit posterior distribution, MCMC provides standard approach to latent variable models.

However, both EM and MCMC methods have prerequisites that are often violated. In practice, not all log-likelihood has finite posterior expectation. Similarly, posterior distributions without tractable behaviour fail the algorithms. Implementing those approaches, at the meantime, could be computationally infeasible, especially when the model structure is complex. In this paper, we establish a unified framework, Maximum Ideal Likelihood Estimation (MILE), remaining valid when traditional methods fail. By treating the latent variables as parameters, naming ``parameterise the latent variables" or ``latent variable parameterisation", MILE estimates the model parameters and latent variables simultaneously, which enjoys consistency, asymptotic distributions, and good empirical performance. Moreover, MILE is asymptotically efficient under regularity. 

Simulation studies indicate that MILE outperforms traditional methods regarding statistical efficiency, computational speed, and prediction accuracy. Overall, MILE, borrowing strength from modern optimisation techniques, is promising and efficient, and serves as an alternative to traditional methods, with novelty and direct interpretability. Note that MILE does not explicitly require a prior distribution, but it could be readily adopted in Bayesian settings for inference. Nevertheless, MILE is compatible with cutting-edge techniques, e.g., Variational inference \citep{Jordan1999, Wainwright2008, VBReview}, and approximate likelihood methods \citep{Lindsay1988, Varin2011, Reid2013, VecchiaReview2021}.

This paper is arranged by properties of MILE. Section \ref{sec:Background} summarises the framework and compares it with traditional methods. Section \ref{sec:ComputationalProcedures} discusses the implementation details, such as numerical algorithms. Section \ref{sec:Inference} presents large sample results of consistency and asymptotic distributions. Simulation studies and illustrative data example are completed in Sections \ref{sec:Simulation}. Section \ref{sec:Discussion} provides the discussion and future works. Additionally, computation algorithms and technical proofs are presented in the Supplement Materials.

\section{Background}
\label{sec:Background}

\subsection{Traditional Approaches to Model with Latent Variables}
The EM algorithm \citep{EMHartley1958, EM1977} is a powerful technique to latent variable models. It shows better stated consistency after \cite{EMWu1983} corrected the flaw in the original proof of the procedure. Denote observed data as $\X \in \mathcal{X}$, latent variables as $\Z \in \mathcal{Z}$ and parameters as $\btheta \in \bTheta$, with joint probability $f(\X, \Z|\btheta)$ and marginal probability $f(\X|\btheta)$. To maximise the marginal likelihood,

$$
L(\btheta; \X) = f(\X|\btheta) = \int_{\mathcal{Z}}f(\X, \Z| \btheta) d\Z,
$$
the EM algorithm, applied to the conditional expectation of log-likelihood followed by a maximisation step, generates sequences of estimators $\{\widehat{\btheta}^{(t)}\}$. The EM estimator, which is the limit of the sequences $\widehat{\btheta} = \lim_{t \rightarrow \infty} \widehat{\btheta}^{(t)}$, is exactly the MLE. Due to the preferred property, the EM algorithm becomes a widely-used technique, but it relies on the prerequisite of a finite conditional expectation of log-likelihood. If the conditions are violated, the EM algorithm does not work because of the failure in the ``E" step. The log-Cauchy Mixture Model, provided in Simulation \ref{Simulation:logCauchy}, is an example whose log-likelihood expectation is infinite, where the EM algorithm is inapplicable.

Other numerical issues also lead to problems. Bayesian Segment Regression (BSR) in Simulation \ref{Simulation:BayesianSegmentRegression} has finite posterior expectation but without analytical expression, rendering its numerical solution computationally intractable, which implies the EM algorithm not suitable for BSR.

Over the past decades, modifications are proposed to overcome computation challenges. The Monte Carlo Expectation Maximisation (MCEM) \citep{MCEM2001} is an important generalization, approximating the expectations via Monte Carlo samples. See \cite{MCEMreview2024} for review. Other alternatives, such as the stochastic EM algorithm by \cite{stEM2000}, simulated EM algorithm by \cite{RUUD1991} and Expectation Conditional Maximisation (ECM) algorithm by \cite{MXL1993}, are remarkable in different scenarios. The Majorise-Maximisation (MM) algorithm \citep{KennethMMA} is another class of algorithms to solutions, but requiring surrogate functions instead of finite expectations, while the function selection could be difficult and limits the flexibility. Undoubtedly, when the original EM algorithm fails due to numerical issues, MCEM and the related methods could be alternatives, but they still underlies other fundamental conditions. In the examples of log-Cauchy Mixture Model and BSR, the alternatives lose feasibility, and thus, it shows that the EM-type algorithms are not applicable in many scenarios.

The limitations of moment-basis also explains the major obstacles to other alternative classes. Under the Bayesian settings, MCMC is a class of methods that is preferred by many researchers. See \cite{GelfandSmith1990} and \cite{MCMChistory} for development and summaries. MCMC aims to approximate the posterior distribution, and many numerical approaches strengthen the performance of MCMC. The acceptance-rejection algorithm \citep{ARMethod}, Metropolis-Hastings algorithm \citep{Metropolis1953,Hastings1970} and Gibbs sampling \citep{GibbsSampling} often shows in MCMC approaches, improving the results of estimators. By MCMC, sampling from the posterior distribution, we get series of dependent estimators, where empirical methods of burning-in (exclude initial part) and thinning (use sub-sequences) are typical to reduce dependence. However, the MCMC techniques are usually slow, taking it longer to derive a stable result \citep{Shen2010}. MCMC also requires parameters with finite posterior expectation, which is a weakness leading to problems including slow sampling, dependent sequences and, similar to the EM algorithm, finite expectation prerequisite.

Briefly, the EM algorithm and MCMC both require finite (posterior or conditional) expectation in some ways. Furthermore, the EM algorithm does not estimate latent variables directly, so alternatives, which includes the latent value estimates, are preferred, if the latent variable enjoys strong importance, such as cases of cluster labels and latent factors for structural equation. To overcome the limitations of traditional methods, a unified framework, with corresponding optimisation algorithms in different scenarios, is raised and introduced in this paper, providing a novel perspective of estimation. We name the framework Maximum Ideal Likelihood Estimator (MILE), as the proposed method maximises the likelihood in ideal cases, by a similar idea of MLE by \cite{fisher1912}. The meaning of ``ideal cases" will be explained later in this section.

\subsection{Motivations}

The motivation originates from a fundamental question that what necessitates the development of the EM algorithm, and at a broader level, the latent variable models. In practice, the central task of statistics is to uncover patterns from observations, so that independence assumptions are typically introduced for tractability and interpretability. The assumptions, though idealised, have been critical to classical methodology and remains a guiding principle.

However, many datasets exhibit strong dependence, where independence assumptions lead to bias. Latent variables offers a natural mechanism to explain the dependence, which maintains interpretability and accounts for their wide use in diverse applications. From a technical perspective, latent variable models provide a route to the MLE, clarifying the necessity of the EM algorithm and its alternatives. Although the MLE is the benchmark for statistical efficiency, the EM algorithm is a computational surrogate throughout with strong limitations. Crucially, if latent variables are observed, the resulting estimators would outperform the marginal MLE, implying the MLE is not the only solution to the problems.

Methodologically, there arises another canonical approach. The two-stage methods \citep{Joe2005,TF2024} divide the parameter estimations into sequential steps. Some quantities are estimated in the first stage, and conditionally on which, the subsequent parameters estimators are derived. However, the sequential construction is prone to error accumulation risks, which undermines efficiency. Hence, these challenges call for alternatives that integrate both stage simultaneously, highlighting the value of other approaches, such as MILE that directly targets the ideal likelihood.

\subsection{MILE Framework Set-up}
MILE expands support of numerical solutions to complicated models. As stated in the latent variable models, we maximise the marginal likelihood

$$
L(\btheta; \X) = f(\X|\btheta) = \int_{\mathcal{Z}}f(\X, \Z| \btheta) d\Z,
$$
which could be viewed as a moment. Note that derivation of $L(\btheta; \X)$ typically involves a challenging integration, which motivates the use of the EM algorithm. However, for the sole purpose of constructing estimators, integration is not the only available approach. Alternatively in MILE, we maximise the integrand $f(\X, \Z| \btheta)$ instead of the integration. Under likelihood format of $L(\Z, \btheta | \X) = f(\X, \Z| \btheta)$, the estimators are

\begin{equation}
\begin{aligned}
    \big(\widehat{\btheta}, \widehat{\Z}\big) & = \argmax_{\Z \in \mathcal{Z}, \btheta \in \bTheta} L(\Z, \btheta | \X) = \argmax_{\btheta \in \bTheta, \Z \in \mathcal{Z}} f(\X, \Z | \btheta)\\
    & = \argmax_{\btheta \in \bTheta, \Z \in \mathcal{Z}} \log f(\X, \Z | \btheta)\\
    &:= \argmax_{\btheta \in \bTheta, \Z \in \mathcal{Z}}  \ell(\Z, \btheta ; \X).
\end{aligned}
\label{TargetLoss}
\end{equation}
Although $L(\Z, \btheta | \X)$ has similar structure of likelihood, note that the exact latent values are not observed. We name $L(\Z, \btheta | \X)$ as the ``ideal" likelihood, because it equals to joint likelihood $L(\Z, \btheta | \X) = L(\btheta | \X, \Z)$ in ``ideal" cases where we know the latent values. The motivation of MILE is to parameterise latent variables $\Z$, i.e. to treat $\Z$ as parameters, whose support $\mathcal{Z}$ plays as the parameter space.

MILE overcomes some limitations of the EM algorithm and MCMC, carrying some superior performance. Finite expectations assumption is no longer a fundamental factor in MILE framework, implying that MILE usually exists even when the traditional methods fail. MILE could be adopted in Bayesian settings, whereas it operates without specification of prior distributions, unlike MCMC which inherently includes them. Remarkably, MILE outperforms the major competitors in terms of computational speed (See Section \ref{sec:Simulation}), achieving identical asymptotic variances to the marginal MLE (See verification in Supplement 4) in some scenarios. Comparison details could be found in Comparison Table in Section \ref{SectionCompTable}. Considering the possible large number of ``parameters" under MILE framework, some numerical algorithms are included and presented, from which MILE incorporates strength.

\subsection{Comparisons between Common Methods}
\label{SectionCompTable}

As previously discussed, prerequisites of traditional methods are not necessary for MILE. MCMC underlies a Bayesian setting of finite posterior expectation, with mechanism to generate random numbers. Similarly, EM-like methods require a finite conditionally expected log-likelihood, which can be calculated (or generated) with computational feasibility. Besides, more conditions might be introduced case-wisely.

MILE offers flexibility and does not rely on expectation-based constrains. MILE remains applicable for estimation and inference in settings where MCMC and the EM algorithm fail. Table \ref{ComparisonTable} provides a comparative overview among MILE, the EM algorithm and MCMC, in terms of implementation and empirical performance.

Column 1 to 4 denote the problem settings: $\pi_{\btheta}$, whether parameters there is prior distribution; $\partial_{\Z}$, whether latent variate (sub-)gradient is finite and numerically feasible; $\Z_{\bpi^*}$, whether latent variables posterior distribution can be generated; $\mathbbm{E}_{\Z_{\bpi^*}}$, whether latent variables posterior expectation is finite and numerically feasible.

For additional conditions and applicability, denote: \ding{51}, conditions meet or algorithm applicable; \ding{53}, conditions not meet either algorithm not applicable; $\forall$, the condition have no impact; ?, conditions and algorithm to be determined. The last column denotes computational speed: $\gtrdot$, faster than competitors; $\lessdot$, slower than competitors; ?, to be determined.

Notice that, as long as there is a density function or massive function, MILE usually exists despite of the speed.

\subsection{Philosophical Ideas and Mathematical Nature}
\label{Section:MeasureTransfer}

It might be noticed that MILE is somehow similar to Hierarchical Likelihood (H-likelihood) by \cite{1996Nelder}, but they are different as a framework.

\cite{1996Nelder} discussed an extended likelihood application in Hierarchical Generalized Linear Model (GLM), whereas the heuristic idea was not completely developed. Its inference strongly relies on the GLM structure, instead of a general distribution. As a result, the reliance to model obscures the mathematical nature of the idea, thereby leading to insufficient justification within a statistical perspective. Furthermore, their statement of selection of target function, between marginal likelihood and H-likelihood, is grounded to an optimisation consideration, relying on reasoning that is narrow in scope and needs deeper theoretical insight. 

This outline might lead to systematic misjudgment. Due to few convincing statistical perspectives, the emphasis on optimisation hinders the identification to proper competitors. \cite{2011Meng} criticised the point in his speech that people understand ``minimizing a loss", but more care should be assigned to why a different loss is selected. Unfortunately, after years of development, there is not enough demonstration in relevant researches.

As a framework, H-likelihood remains incomplete and requires refinement. Target functions could  exhibit various mathematical property, while H-likelihood fails to cover many of them due to methodological limitations, which in part explains the narrow scope of its application and disclosure. In review of \cite{Lee2021}, H-likelihood remains highly confined to GLM and its extension, under strong continuity restrictions, such as in hierarchical survival models. Regardless of the scope, the idea potential deserves further exploration.

Furthermore, ``h" for hierarchical is not a proper name to represent the philosophy of the framework. In non-hierarchical models with complicated structure, development could be pursued by extending the idea. Thus, MILE approaches the problem in a more general prospective, not relying on specific structural assumption. Restricted to (H-)GLMs, H-likelihood could be viewed as a special case of MILE.

``Latent variable parameterisation" constitutes the most crucial procedure that is intrinsically connected to its underlying philosophy. Parameterising latent variables by treating $\Z$ as unknown constants instead of random variables, should not be viewed solely a numerical strategy of optimisation, as statistical interpretations is still required.

Parameterisation is fundamentally a probability measure transfer. Denote the original probability system defined on $(\Omega, \mathcal{B}, \mathbbm{P})$, where $\Omega = \mathcal{X}$, $\mathcal{B}$ is the $\sigma$-algebra of $\Omega$ on some probability function $\mathbbm{P}$. Notice that $\Z$ is unobserved, rendering the $\sigma$-algebra unmeasurable generated by $(\mathcal{X}\times\mathcal{Z})$. After parameterisation, $(\Omega, \mathcal{B}, \mathbbm{P})$ is transferred to a new measure $(\Omega, \mathcal{B}_0, \mathbbm{P}_0)$. By the Radon-Nikodym Theorem (Theorem 6.10 of \cite{rudin86real}), we transfer any $B \in \mathcal{B}$ to some $B_0 \in \mathcal{B}_0$ by

$$
B  = \int_{B} d\omega= \int_{B_0}\frac{d \omega}{d \omega_0} d\omega_0.
$$

Reversely, for any $A_0 \in \mathcal{B}_0$, its original event $A$ can be expressed as
$$
A_0 = \int_{A_0} d\omega_0= \int_{A}\frac{d \omega_0}{d \omega} d\omega.
$$

When latent variables are treated as constants, the analysis is assumed that the latent variables are fixed whichever known or unknown. Equivalently, it corresponds to work under a conditional probability measure given $\Z$. This perspective is analogous to the model misspecifcation scenarios as outlined below. Suppose there are identically and independently distributed (i.i.d) samples from some distribution with pdf $f(\theta_1, \theta_2)$. When the value of $\theta_2$ is wrongly given as $\theta_2 = \theta_2'$, we could still get estimates of $\widehat{\theta}_1$. The inferences could be derived under structure of Godambe's Sandwiches \citep{Godambe1960}, and \cite{Joe2005} analysed one of its application.

Accordingly, the Radon-Nikodym derivative follows the expression as $d \omega_0/d \omega = f(\X| \Z, \btheta)/f(\X, \Z| \btheta)$ or $d \omega/d \omega_0 = f(\X, \Z| \btheta)/f(\X| \Z, \btheta)$. Suppose there is an estimator $\widehat{\btheta}(\X)$ of $\btheta$, and for any event $A \in \mathcal{B} \cap \mathcal{B}_0$, we have probability as
\begin{equation}
    P_{\widehat{\btheta}(\X) | \btheta}(A)=P_{\widehat{\btheta}(\X) | \btheta}\left(\int_A d \omega\right)=\int_A d F(\X | \btheta)=\int_A \int_{\mathcal{Z}} d\Z d F(\X, \Z | \btheta).
    \label{OriginalProb}
\end{equation}

The corresponding transferred probability is
\begin{equation}
    \begin{aligned}
        P_{\widehat{\btheta}(\X) | \btheta, \Z}(A) & = P_{\widehat{\btheta}(\X) | \btheta}(A_0) = P_{\widehat{\btheta}(\X) | \btheta}\left(\int_A \frac{d \omega_0}{d \omega} d \omega\right)\\
        &=\int_A \int_{\mathcal{Z}} \frac{d \omega_0}{d \omega} d \Z d F(\X, \Z | \btheta) =\int_A \int_{\mathcal{Z}} \frac{d \omega_0}{d \omega} f(\X, \Z | \btheta) d\Z d\X \\
        & =\int_{A} \int_{\mathcal{Z}} f(\X | \btheta, \Z) d \Z d \X=\int_{A} \int_{\mathcal{Z}}  d \Z d F(\X | \btheta, \Z)
    \end{aligned}
    \label{TransferredProb}
\end{equation}
\eqref{OriginalProb} and \eqref{TransferredProb} show how event probability differ before and after measure transfer. We will discuss the estimator distribution given true value of $\Z$ and $\btheta$, i.e., $(\widehat{\btheta}, \widehat{\Z})|\{\Z, \btheta, \X\}$,  in the following sections. If the marginal distribution is primarily interested, it can be obtained through a slight backward transfer via Bayes' formula.

To establish a conprehensive framework, we expand applications of MILE in scenarios of categorical and/or non-differentiable latent variables, including clustering and related challenging problems. Section \ref{sec:ComputationalProcedures} presents corresponding algorithms, rendering the framework practical in a wide range of numerical issues.

\section{Implementation and Methodology}

\label{sec:ComputationalProcedures}

\subsection{Proposed Methodology}

Usually given $\Z$ and $\X$, it isn't hard to obtain the estimator $\widehat{\btheta} | (\Z, \X)$. Thus we replace the latent random variables $\Z$ by $\widehat{\Z}$, conducting the shift onto $\widehat{\btheta} | (\widehat{\Z}, \X)$.

Furthermore, we parameterise $\Z$ to solve (\ref{TargetLoss}). If $\Z$ has partial gradients,

\begin{equation}
    \frac{\partial \ell(\Z, \btheta ; \X)}{\partial \Z} = \frac{\partial \big\{\log f(\X | \Z, \btheta) + \log f(\Z | \btheta)\big\}}{\partial \Z} = \frac{\partial \log f(\X | \Z, \btheta)}{\partial \Z} + \frac{\partial \log f(\Z | \btheta)}{\partial \Z};
    \label{ZGradient}
\end{equation}

\begin{equation}
    \frac{\partial \ell(\Z, \btheta ; \X)}{\partial \btheta} = \frac{\partial \big\{\log f(\X | \Z, \btheta) + \log f(\Z | \btheta)\big\}}{\partial \btheta} = \frac{\partial \log f(\X | \Z, \btheta)}{\partial \btheta} + \frac{\partial \log f(\Z | \btheta)}{\partial \btheta}.
    \label{thetaGradient}
\end{equation}
Estimators of $\Z$ and $\btheta$ can be derived by solving (\ref{ZGradient}) and (\ref{thetaGradient}), while Algorithm 2 is employed when no close form is available. Note that, nevertheless, neither partial gradients is necessarily required, as without which there still are maximisers.

In practice, obtaining $\widehat{\Z}$ could be challenging, because of, for example, absense of (sub-)gradients. Hence, derivative-free scenarios are of our particular interest. Due to the complicated expression in \eqref{TargetLoss}, different methods are required, with empirical details presented in Section \ref{sec:Simulation}. For now, corresponding methods, serving to derive (global or local) maximisers of the ideal likelihood function, thus adequately cover scenarios which are frequently encountered in practice.

\subsection{Computational Procedures}
As assumed, $\btheta$ retains favourable properties under idea likelihood functions, we focus on scenario of latent variables $\Z$.

\subsubsection{Continuous Latent optimisation}

For $\Z$ with continuous support, two cases represent: differentiable log ideal likelihood with respect to $\Z$; and non-differentiable log ideal likelihood with respect to $\Z$.

When the log ideal likelihood is differentiable, the optimisation problem can be addressed by Block Coordinate Ascending (BCA) in Algorithm 2. The initials of $\btheta$, required by BCA, are denoted as $\btheta^{(0)}$, after which the maximisers of \eqref{BCA} are computed iteratively until convergence.

\begin{equation}
    \begin{aligned}
        \Z^{(t + 1)} &= \argmax_{\Z \in \mathcal{Z}} \ell(\btheta^{(t)}, \Z; \X)\\
        \btheta^{(t + 1)} &= \argmax_{\btheta \in \Theta} \ell(\btheta, \Z^{(t+1)}; \X)
    \end{aligned}
    \label{BCA}
\end{equation}

When the log ideal likelihood is not differentiable, we employ the Genetic Algorithm (GA) to estimate latent values. GA is a special class of evolutionary algorithm which does not rely on gradient. It is well-suited for complex settings such as change point detection \citep{Davis2006}. See \cite{Holland1992} and \cite{Eiben2015} for reviews. We further propose a hybrid GA in Algorithm 1 which estimates latent variables and parameters simultaneously, where chromosomes are decoded as $\widehat{\Z}$.  Under the chromosome interpretation, \eqref{thetaGradient} is reformulated to \eqref{thetaGradientOnZhat} for the zero points as parameter estimators,
\begin{equation}
    \frac{\partial \ell(\btheta; \widehat{\Z}, \X)}{\partial \btheta} = \frac{\partial \big\{\log f(\X | \widehat{\Z}, \btheta) + \log f(\widehat{\Z} | \btheta)\big\}}{\partial \btheta} = \frac{\partial \log f(\X | \widehat{\Z}, \btheta)}{\partial \btheta} + \frac{\partial \log f(\widehat{\Z} | \btheta)}{\partial \btheta}.
    \label{thetaGradientOnZhat}
\end{equation}
Both ideal likelihood $L(\widehat{\btheta}, \widehat{\Z}; \X)$ or log ideal likelihood $\ell(\widehat{\btheta}, \widehat{\Z}; \X)$ can be applied as the fitness function. For numerical stability, logarithm forms are generally preferred.

Although the algorithm is organized in a step-wise manner, both $\widehat{\Z}$ and $\widehat{\btheta}$ are simultaneously obtained through an entire evaluation of fitness function. To improve the empirical performance of GA, initialization strategies such as Voronoi Partition \citep{Shimosaka} might be applied. In high-dimensional latent variables settings, the improvement is particularly crucial to use hybrid GA. Additional numerical techniques are acknowledged but left for future researches.

\subsubsection{Categorical Optimisation}
To optimise log ideal likelihood with categorical latent $\Z$, we introduce the following definitions.

\begin{definition}{\textbf{(Slice-wise Convex Function)}}
    Suppose a real-number function $f(\X, \Y)$. Domain of $\X = (X_1, X_2, \cdots, X_n)$, $\mathcal{D}(\X)$, is categorical, and domain of $\Y = (Y_1, Y_2, \cdots, Y_m)$, $\mathcal{D}(\Y) \subset \mathbbm{R}^m$, is a subset of real region, where
    $$
    \mathcal{D}(\X) = \mathcal{D}(X_1) \times \mathcal{D}(X_2) \times \cdots \times \mathcal{D}(X_n)
    $$
    for all $i \in \{1, 2, \cdots, n\}$, $X_i \in \mathcal{D}(X_i) = \{x_{i, 1}, x_{i, 2}, \cdots, x_{i, n_i}\}$. $f(\X, \Y)$ is \textbf{slice-wise convex}, if $f(\x, \Y)$ is convex, $\forall$ $\x \in \mathcal{D}(\X)$.
    \label{SlicewiseConvexFunction}
\end{definition}

Definition \ref{SlicewiseConvexFunction} characterise the convexity of functions with categorical inputs. Note that the log ideal likelihood function of Gaussian Mixture Model (GMM) is slice-wise convex in particular. It is natural to assume that the log ideal likelihood with categorical latent variables and continuous parameters satisfies slice-wise convexity. Similar assumptions have been adopted discrete parameters models \citep{Choirat2012,Ma2023}.

Furthermore, we need a definition of neighbourhood relationship.

\begin{definition}{(Categorical Optimisational Neighbourhood)}
    Suppose categorical vectors, $\x = (x_1, x_2, \cdots, x_n)$ and $\y = (y_1, y_2, \cdots, y_n)$, share the same domain $\mathcal{D}$.
    
    $\x$ and $\y$ are \textbf{optimisational neighbours}, notated as $\x \sim \y$, if $\exists$ unique $j \in \{1, 2, \cdots, n\}$, $x_j \neq y_j$, and $\forall$ $i \neq j$, $x_i = y_i$.
    \label{CategoricaloptimisationalNeighbourhood}
\end{definition}

By Definitions \ref{SlicewiseConvexFunction} and \ref{CategoricaloptimisationalNeighbourhood}, we developped an optimisation algorithm for categorical inputs, termed Stepwise Categorical Progress (SCP) as outlined in Algorithm \ref{SCP}.

\begin{algorithm}[hbtp]
\caption{Stepwise Categorical Progress}\label{SCP}
\begin{algorithmic}[1]
\Require Observations, Slice-wise convex fitness $\ell(\Z, \btheta ; \X)$, MaxIter = $m$
\State \textbf{Initialise} count = 1, latent variables $\Z$
\While{count < $m$}
\State \textbf{Initialise} opt = 0
\For{$i$ \textbf{{\color{blue}in}} $1:Population$}
\For{$k$ \textbf{{\color{blue}in}} $1:n_i$}
\State \textbf{Set} a temporary latent vector $\Z^\prime = \Z$, and \textbf{Set} $\Z^\prime(i) = k$
\State Calculate $\widehat{\btheta} = \underset{\btheta \in \bTheta}{\argmax} \ell(\Z, \btheta ; \X)$, $\widehat{\btheta^\prime} = \underset{\btheta \in \bTheta}{\argmax} \ell(\Z^\prime, \btheta; \X)$
\State Calculate Old-Fitness $\ell(\Z, \widehat{\btheta}; \X)$ and New-Fitness $\ell(\Z^\prime, \widehat{\btheta^\prime} ; \X)$

\If{$\ell(\Z^\prime, \widehat{\btheta^\prime} ; \X)$ > $\ell(\Z, \widehat{\btheta}; \X)$}
\State \textbf{Replace} $\Z(i) = k$
\State opt = 1
\State \textbf{Next} $k$
\EndIf

\EndFor
\EndFor
\If{opt == 0}
\State {\color{blue} \textbf{Break}} \textbf{While}
\EndIf
\State count = count + 1
\EndWhile
\State \Return $\Z$
\end{algorithmic}
\end{algorithm}

SCP optimises the target function in two stages by: searching the neighbours of a candidate solution for improvement; traverse all solutions until no neighbour yields higher fitness. Simulation \ref{Simulation:GMM} illlustrates the reliability of SCP. Definition \ref{LocalMaximumOfCategoricalFunction} formalises the local maximiser of a categorical function. Because, obviously, categorical function with finite domain always admits at least one local maximum, SCP reaches a local maxima.

\begin{definition}{(Local Maximum of Categorical Function)}
    For a function $f(\X)$ with categorical domain $\mathcal{D}$, $\x \in \mathcal{D}$ is a \textbf{local maximiser}, if 
    $$
    f(\x) \geq f(\y)
    $$
    for $\forall$ $\y \in \mathcal{D}$ satisfying $\x \sim \y$.
    \label{LocalMaximumOfCategoricalFunction}
\end{definition}

\begin{theorem}
    SCP converges to $\widehat{\Z}$, the unique local maximiser of the fitness function $g(\Z)$, if $g(\Z)$ is slice-wise convex with categorical input and finite domain.
    \label{SCPConvergence}
\end{theorem}

\begin{proof}

Let initial of SCP be $\Z$, and denote the optimiser after steps to be $\Z_1, \Z_2, \cdots \Z_n \cdots$, with eventual optimiser $\widetilde{\Z}$. Suppose that $\widetilde{\Z}$ is not a local maximiser, so there exists another $\Z_{next} \sim \widetilde{\Z}$, such that
$$
    g(\Z_{next}) > g(\widetilde{\Z}).
$$
        
In this case, SCP doesn't terminate at $\widetilde{\Z}$, contradicting to the assumption of the final output. Furthermore, if $\widetilde{\Z} \neq \widehat{\Z}$, $g(\Z)$ admits at least two distinct local maxima, contradicting to uniqueness assumption. Thus, SCP converges to $\widehat{\Z}$.
\end{proof}

Theorem \ref{SCPConvergence} presents a sufficient condition under which SCP admits the global maximiser. For easier notations in the following sections, denote $\widehat{\btheta}(\Z) = \underset{\btheta \in \bTheta}{\argmax} \ell(\Z, \btheta ; \X)$ and set $g(\Z) = \ell\big(\Z, \widehat{\btheta}(\Z) ; \X \big).$

\section{Inference with Interpretation}

\label{sec:Inference}

\subsection{Assumptions}
By parameterising the latent variables $\Z$, we obtain the estimators $\widehat{\Z}$ requiring a perspective for inference. It ought to proceed as if $\Z$ are constants, instead of random variables. Although classical inference techniques are still applicable, they should be conditionally based on $\Z$. Formally, it requires a transfer from an unconditional measure to a conditional probability measure $(\widehat{\btheta}, \widehat{\Z})|\Z(\omega)$, for some random events $\omega \in \Omega$. For notational convenience, abbreviate it as $(\widehat{\btheta}, \widehat{\Z})|\Z$ hereafter.

Regularity, the fundamental and widely-used assumption, is required for inference, based on which, we conduct the theoretical analysis for MILE. Because MILE closely relies on specific probability measures in Section \ref{Section:MeasureTransfer}, distinct regularity is introduced. Assumption \ref{ConditionalRegularity} concerns the conditional measure $\X | \Z, \widehat{\Z}$, while Assumption \ref{MarginalRegularity} concerns the marginal measure $\X | \Z$. Regardless of the assumption item amount, they are natural and generally satisfied in empirical studies.

\begin{assumption}\textbf{(Regularity of Ideal Likelihood Function)}
    In latent variable model, let $\X \in \mathcal{X} \subset \mathbbm{R}^n$ be a random vector, where parameter space $\bTheta$ and latent variable space $\mathcal{Z}$ are subsets of real region. The conditional probability function $f(\X | \Z, \btheta): \mathcal{X} \times \bTheta \times \mathcal{Z} \rightarrow \mathbbm{R}$ and ideal likelihood function $h(\X, \Z | \btheta): \mathcal{X} \times \bTheta \times \mathcal{Z} \rightarrow \mathbbm{R}$ are real-value functions. With some $N$, assume:

    \begin{enumerate}[label=(\alph*)~]
        \item \textbf{(Compactness)} $\bTheta \times \mathcal{Z}$ is compact; \label{assumption_compact}
        
        \item \textbf{(Positive Value)} for $\forall$ $(\X, \btheta, \Z) \in \mathcal{X} \times \bTheta \times \mathcal{Z}$, $h(\X, \Z | \btheta) > 0$, $f(\X | \Z, \btheta) > 0$;

        \item \textbf{(Measurability)} for $\forall$ $(\btheta, \Z) \in \bTheta \times \mathcal{Z}$, $f(\X | \Z, \btheta)$, $h(\X, \Z | \btheta)$ are $\mathcal{X}$-measurable;\label{assumption_measure}

        \item \textbf{(Continuity)} for $\forall$ $\X \in \mathcal{X}, \Z \in \mathcal{Z}$, $h(\X, \cdot, \Z)$ is continuous on $\bTheta$;\label{assumption_continuous}
        
        \item \textbf{(Kullback–Leibler (K-L) Divergence)} for $\forall$ $(\btheta, \Z), (\widehat{\btheta}, \widehat{\Z})\in \bTheta \times \mathcal{Z}$, $KL(\Z, \btheta, \widehat{\Z}, \widehat{\btheta})=\frac{1}{N}\int_\mathcal{X} f(\X | \Z, \btheta) \log \frac{f(\X | \Z, \btheta)}{h(\X, \widehat{\Z} | \widehat{\btheta})} d \X<\infty$; \label{assumption_KLDiversity}

        \item \textbf{(Fubini's Interchange)} for $\forall$ $\Z \in \mathcal{Z}$ and $\widehat{\Z} \in \mathcal{Z}$, $
            \frac{\partial}{\partial \btheta} \int f(\X | \Z, \btheta^c) d \X=\int \frac{\partial}{\partial \btheta} f(\X | \Z, \btheta^c) d \X$ and $\frac{\partial}{\partial \btheta} \int h(\X, \widehat{\Z} | \btheta) d \X=\int \frac{\partial}{\partial \btheta} h(\X, \widehat{\Z} | \btheta) d \X$;
        \label{assumption_Fubini}
    \end{enumerate}
    \label{FunctionalRegularity}
\end{assumption}

\begin{assumption}\textbf{(Conditional Regularity)}
   Under conditional measure $\X|\widehat{\Z}, \Z$ and by notations and items in Assumption \ref{FunctionalRegularity}, further assume:

    \begin{enumerate}[label=(\alph*)~]
        \item \textbf{(Score Variance)} for $\forall$ $(\btheta^c, \Z), (\btheta, \widehat{\Z})\in \bTheta \times \mathcal{Z}$, $||\V(\btheta)||_{\infty} < \infty$ where
        $$\V(\btheta) = \int_\mathcal{X} \frac{f(\X | \Z, \btheta^c)}{N}\left(\frac{\partial \log h(\X,  \widehat{\Z} | \btheta)}{\partial \btheta}\right)\left(\frac{\partial \log h(\X,  \widehat{\Z} | \btheta)}{\partial \btheta}\right)^\top d \X;$$
        \label{assumption_ScoreVariance}

        \item \textbf{(Curvature)} for $\forall$ $(\btheta^c, \Z), (\btheta, \widehat{\Z})\in \bTheta \times \mathcal{Z}$, 
        $$\A(\btheta) = \frac{1}{N} \int_\mathcal{X} f(\X | \Z, \btheta^c)\frac{\partial^2 \log h(\X, \widehat{\Z} | \btheta)}{\partial \btheta^\top \partial \btheta}  d \X,$$
        and $||\A(\btheta)||_{\infty} < \infty$, $||\frac{1}{N}\frac{\partial^2 \log h(\X, \widehat{\Z} | \btheta)}{\partial \btheta\top \partial \btheta}|_{\btheta = \btheta^c}-\A\left(\btheta^c\right)||_{\infty}=\lowero_ p(1);$
        \label{assumption_Curvature}
        
        \item \textbf{(Bounded Skewness Tensor)} $\T(\btheta)$ is the Skewness Tensor of $\log h(\X, \widehat{\Z} | \btheta)$ to $\btheta$, where $T_{i,j,k}(\btheta) = \frac{\partial^3\log h(\X, \widehat{\Z} | \btheta)}{\partial \theta_i \partial \theta_j \partial \theta_k}$, and $\sup\left|\frac{\T(\btheta)}{N} \right|<\infty$; \label{assumption_BST}

        \item  \textbf{(Uniqueness)} for $\forall$ $\X \in \mathcal{\X}$, $\widehat{\Z }\in \mathcal{Z}$, $\log h(\X, \widehat{\Z} | \widehat{\btheta})$ has unique zero point $\widehat{\btheta} \in \bTheta$;
        \label{assumption_UniqueSolution}
        
        \item \textbf{(Asymptotic Normal Score)} for $\forall$ $(\btheta, \widehat{\Z})\in \bTheta \times \mathcal{Z}$, under the conditional measure, $\frac{\partial \log h(\X,  \widehat{\Z} | \btheta)}{\partial \btheta}$ is asymptotic normal.        
        \label{assumption_NormalScore}
    \end{enumerate}
    \label{ConditionalRegularity}
\end{assumption}

\begin{assumption}\textbf{(Marginal Regularity)}
    Under conditional measure $\X|\Z$ and by notations and items in Assumption \ref{FunctionalRegularity}, further denote $\bzeta = (\btheta, \Z)^\top$ and assume:

    \begin{enumerate}[label=(\alph*)~]
        \item \textbf{(Score Variance)} for $\forall$ $(\btheta^c, \Z_0), (\btheta, \Z)\in \bTheta \times \mathcal{Z}$, $||\V(\bzeta)||_{\infty} < \infty$ where
        $$\V(\bzeta) = \int_\mathcal{X} \frac{f(\X | \Z_0, \btheta^c)}{N}\left(\frac{\partial \log h(\X,  \Z | \btheta)}{\partial \bzeta}\right)\left(\frac{\partial \log h(\X,  \Z | \btheta)}{\partial \bzeta}\right)^\top d \X;$$
        \label{assumption_ScoreVariance}

        \item \textbf{(Curvature)} for $\forall$ $(\btheta^c, \Z_0), (\btheta, \Z)\in \bTheta \times \mathcal{Z}$, 
        $$\A(\bzeta) = \frac{1}{N} \int_\mathcal{X} f(\X | \Z_0, \btheta^c)\frac{\partial^2 \log h(\X, \Z | \btheta)}{\partial \bzeta^\top \partial \bzeta}  d \X,$$
        and $||\A(\bzeta)||_{\infty} < \infty$, $||\frac{1}{N}\frac{\partial^2 \log h(\X, \Z | \btheta)}{\partial \bzeta\top \partial \bzeta}-\A\left(\bzeta\right)||_{\infty}=\lowero_ p(1);$
        \label{assumption_Curvature}
        
        \item \textbf{(Bounded Skewness Tensor)} $\T(\bzeta)$ is the Skewness Tensor of $\log h(\X, \Z | \btheta)$ to $\btheta$, where $T_{i,j,k}(\btheta) = \frac{\partial^3\log h(\X, \Z | \btheta)}{\partial \zeta_i \partial \zeta_j \partial \zeta_k}$, and $\sup\left|\frac{\T(\bzeta)}{N} \right|<\infty$; \label{assumption_BST}

        \item  \textbf{(Uniqueness)} for $\forall$ $\X \in \mathcal{\X}$, $\log h(\X, \widehat{\Z} | \widehat{\btheta})$ has unique zero point $(\widehat{\btheta}, \widehat{\Z})\in \bTheta \times \mathcal{Z}$; 
        \label{assumption_UniqueSolution}
        
        \item \textbf{(Asymptotic Normal Score)} for $\forall$ $(\btheta, \Z) \in \bTheta \times \mathcal{Z}$, under the marginal measure $\X|\Z$, $\frac{\partial \log h(\X,  \Z | \btheta)}{\partial \bzeta}$ is asymptotic normal.        
        \label{assumption_NormalScore}
    \end{enumerate}
    \label{MarginalRegularity}
\end{assumption}

Assumption \ref{FunctionalRegularity} specifies the required fundamental property of the probability functions, complemented by Assumption \ref{ConditionalRegularity} and Assumption \ref{MarginalRegularity}, which impose the conditional and marginal regularity respectively. Although item \ref{assumption_NormalScore} seems strong, it is widely supported by empirical and theoretical results. In i.i.d samples, the Central Limit Theorem (CLT) applies. Furthermore, for weakly dependent data, results such as mixing CLT \citep{Jenish2009}, dependency graph CLT \citep{janson1988MRF} and Bernstein–von Mises Theorem \citep{Lecam1986} for Bayesian inference, could be invoked. Item \ref{assumption_NormalScore}to in necessary for asymptotic normality, but without which non-normal asymptotic distributions remains possible.

$N$, in Assumption \ref{FunctionalRegularity}, is critical in the assumptions. In i.i.d settings, $N$ coincides with the sample size. Under dependence, however, the equivalence generally fails. Broadly, $N$ could be regarded as a function of the sample size, and in many scenarios, $N$ exhibits a power growth, formalized in Assumption \ref{assumption_CovergeRate}.

\begin{assumption}\textbf{(Estimable Convergence Rate)}
    Let $\X \in \mathcal{X} \subset \mathbbm{R}^n$ be a random vector, parameter space $\bTheta$ and latent variable space $\mathcal{Z}$ be subsets of real region. $h: \mathcal{X} \times \bTheta \times \mathcal{Z} \rightarrow \mathbbm{R}$ is a real-value function. Under Assumption \ref{ConditionalRegularity} or Assumption \ref{MarginalRegularity}, for $\forall$ $\X \in \mathcal{X}$ and some $k > 1$, 
        $$
        \sup _{(\btheta, \Z) \in \bTheta \times \mathcal{Z}} \left|\frac{1}{N} h(\X, \btheta, \Z) - \frac{1}{N} h(\X, \widehat{\btheta}, \widehat{\Z})\right|  = \lowero_p\left(n^{\frac{1}{k}-1}\right),
        $$
        where $(\widehat{\btheta}, \widehat{\Z}) = \argmax_{(\btheta, \Z) \in \bTheta \times \mathcal{Z}}h(\X, \btheta, \Z)$, and $n$ is the sample size.\label{assumption_CovergeRate}

\end{assumption}

Convergence rate of latent variable models is case-wise specific. In i.i.d scenarios, the rate typically satisfies $k = 2$. Under dependence, smaller $k$ may arise. See \cite{Zhang2004,zhang2006} and \cite{Davis2013} for examples. Generally, convergence rate could be characterized by CLT adapted to dependence structures \citep{PNASWU2005,Jenish2009,Jenish2012}.

MILE may converge to pseudo trues rather than the true parameter values. Assumption \ref{Negligibility} presents conditions of mis-specification. Notably, in dependent data settings, the scaling term $M$ also follows the power growth in Assumption \ref{assumption_CovergeRate}.

\begin{assumption}\textbf{(Negligibility)}
Based on the Assumption \ref{ConditionalRegularity} or Assumption \ref{MarginalRegularity}, assume:
\begin{enumerate}[label=(\alph*)~]
\item \textbf{(Conditional Score Consistency)} for $\forall$ $(\btheta, \Z) \in \bTheta \times \mathcal{Z}$, there exist some $C_1 \in \mathbbm{R}$, such that $\frac{1}{N}\log f(\X | \Z, \btheta)=C_1+\lowero_p(1)$;
\label{assumption_ConditionalScoreConsistency}
\item \textbf{(Negligible Latent Margin)} for $\forall$ $(\btheta, \Z) \in \bTheta \times \mathcal{Z}$, there exist some $C_2 \in \mathbbm{R}$ and $M$, such that $\frac{1}{M}\log f(\Z | \btheta)=C_2+\lowero(1)$, where $M/N \rightarrow 0$;\label{assumption_NegligibleLatentMargin}

\end{enumerate}
\label{Negligibility}
\end{assumption}

\begin{comment}
    \begin{theorem}
    If log ideal likelihood $\log \big(f(\X |\btheta, \Z) f(\Z| \btheta)\big)$ exists and the Weak Law of Large Number holds, its unique maximiser $\widehat{\btheta}$ is consistent to true value $\btheta_0$. That is
    $$
    \widehat{\btheta} \pconverge \btheta_c^*,
    $$
    as $N\rightarrow \infty$, and
    $$\btheta_c^* = \btheta_0 + \O(...).$$
    \label{pconvergence}

    {\color{red}\fbox{We need \textbf{U}WLLN as Lemma 1, direct prove by technical conditions or directly assume it holds.}}
    
    {\color{red} \fbox{Note that we cannot assume $\Z$ on (random) compact set ($\btheta$ is compact).}}
    
    {\color{red} \fbox{May simply assume the ``boundness" of $\log ff$}}
\end{theorem}
\end{comment}

\subsection{Robust Parametric Normality}

Numerical approaches does not always achieve the global maximiser. In particular, the hybrid GA estimator $(\widehat{\btheta}, \widehat{\Z})$ may differ from $\argmax_{(\btheta, \Z) \in \bTheta \times \mathcal{Z}}f(\X, \btheta, \Z)$, although $\widehat{\btheta} = \argmax_{\btheta \in \bTheta }f(\X, \btheta, \widehat{\Z})$ for fixed $\widehat{\Z}$. Accordingly, inference is naturally formulated in terms of the conditional estimator $\widehat{\btheta}|(\widehat{\Z}, \Z)$.

\begin{theorem}{(Robust Conditional Asymptotic Normality)}

    Suppose log ideal likelihood $h(\X, \Z | \btheta)$ and conditional probability function $f(\X |\Z, \btheta)$ satisfies Assumptions \ref{ConditionalRegularity} with some $N$. The asymptotic conditional distribution of $\widehat{\btheta}| (\widehat{\Z}, \Z)$ is
    $$
   \sqrt{N}(\widehat{\btheta} - \btheta^c)| (\widehat{\Z}, \Z) \dconverge N(\textbf{0}, \I_{\textit{c}}^{-1}),
    $$
    where $\btheta^c$ is a pseudo true of $\btheta$, and $\I_{\textit{c}}$ is positive definite.

    \label{ConditionalAsympDist}
\end{theorem}
Inference is conducted under conditional measures, where $\btheta_c$ and $\I_{\textit{c}}$ are constant, while they may appear rand marginallyom. Formally, one may express $\btheta_c:= \btheta_c(\Z, \widehat{\Z})$, implying $\btheta_c(\Z, \widehat{\Z}):= \btheta_c(\Z, \widehat{\Z})|(\Z, \widehat{\Z})$ is constant. This reasoning extends to all subsequent inference. We emphasize again that $N$ may not be the sample size.

Theorem \ref{ConditionalAsympDist} establishes the conditional distribution is robustly asymptotic normal, centered at pseudo true. Theorem \ref{ConditionalAsympDist} legitimises the inference of Simulation \ref{Simulation:BayesianSegmentRegression}, where the information term $\I_{\textit{c}}$ could be expressed as Godambe Matrix \citep{Godambe1960}. Empirically, Monte Carlo procedure is useful to generate the marginal distribution.

Under stronger assumptions on latent variables estimators $\widehat{\Z}$, more refined expressions can be derived.

\subsection{Joint Asymptotic Normality}
$(\widehat{\btheta}, \widehat{\Z})$ could be derived via \eqref{MILEEstDerivation}. Let $\log f(\X | \btheta, \Z) = f_1(\btheta, \Z)$, $\log f\left(\Z | \btheta\right) = f_2(\btheta, \Z)$, and hence $h(\X, \btheta, \Z) = f_1(\btheta, \Z) + f_2(\btheta, \Z)$. 

\begin{equation}
    \begin{aligned}
        (\widehat{\btheta}, \widehat{\Z})&=\argmax_{(\btheta, \Z) \in \bTheta \times \mathcal{Z}} f(\X, \Z | \btheta)=\argmax_{(\btheta, \Z) \in \bTheta \times \mathcal{Z}} L(\btheta; \X, \Z) \\
        & =\argmax_{(\btheta, \Z) \in \bTheta \times \mathcal{Z}} \ell(\btheta, \Z; \X) =\argmax_{(\btheta, \Z) \in \bTheta \times \mathcal{Z}} \big[\log f(\X | \btheta, \Z)+\log f\left(\Z | \btheta\right)\big].
    \end{aligned}
    \label{MILEEstDerivation}
\end{equation}

\begin{theorem}{(Conditional Joint Asymptotic Normality)}
        Suppose log ideal likelihood $h(\X, \Z|\btheta)$ and conditional probability function $f(\X |\Z, \btheta)$ satisfy Assumptions \ref{MarginalRegularity} and Assumption \ref{Negligibility} with some $N$. The asymptotic distribution is
        \begin{equation*}
            \sqrt{N}\left(\widehat{\btheta} - \btheta, \widehat{\Z} - \Z\right) | \Z \dconverge N\left(\textbf{0}, \I^{-1}(\btheta, \Z)\right).
        \end{equation*}    
\label{JointAsymDist}
\end{theorem}
\begin{comment}
    Notice that 
\begin{equation}
            \I(\btheta, \Z)=\frac{1}{N} \left(\begin{array}{cc}
        -\E \big[ \frac{\partial^2 f_1(\btheta, \Z)}{\partial \btheta^\top \partial \btheta} - \frac{\partial^2 f_2(\btheta, \Z)}{\partial \btheta^\top \partial \btheta} \big] & -\E \big[ \frac{\partial^2 f_1(\btheta, \Z)}{\partial \btheta^\top \partial \Z} -\frac{\partial^2 f_2(\btheta, \Z)}{\partial \btheta^\top \partial \Z} \big] \\
        -\E \big[ \frac{\partial^2 f_1(\btheta, \Z)}{\partial \btheta^\top \partial \Z} - \frac{\partial^2 f_2(\btheta, \Z)}{\partial \btheta^\top \partial \Z} \big] & -\E \big[ \frac{\partial^2 f_1(\btheta, \Z)}{\partial \Z^\top \partial \Z} - \frac{\partial^2 f_2(\btheta, \Z)}{\partial \Z^\top \partial \Z} \big]
        \end{array}\right).
        \label{I_theta0_Z}
        \end{equation}
\end{comment}

Proof of Theorem \ref{JointAsymDist} follows standard asymptotic arguments \citep{Shao2003,Miller2021}, with key details provided in Supplement 3. Assumption \ref{Negligibility} is not required for the asymptotic normality; in its absence, the estimator converges to a pseudo true. In terms of the information, the blocked information matrix format could be used
$$
\I(\btheta, \Z)= \left(\begin{array}{cc}
\I_{11} & \I_{12} \\
\I_{12}^{\top} & \I_{22}
\end{array}\right),
$$

\begin{comment}
\begin{eqnarray*}
\I(\btheta, \Z)&=& \left(\begin{array}{cc}
\I_{11} & \I_{12} \\
\I_{12}^{\top} & \I_{22}
\end{array}\right)\\
&=&
\left(\begin{array}{cc}
        -\E \big[ \frac{\partial^2 f_1(\btheta, \Z)}{\partial \btheta^\top \partial \btheta} - \frac{\partial^2 f_2(\btheta, \Z)}{\partial \btheta^\top \partial \btheta} \big] & -\E \big[ \frac{\partial^2 f_1(\btheta, \Z)}{\partial \btheta^\top \partial \Z} -\frac{\partial^2 f_2(\btheta, \Z)}{\partial \btheta^\top \partial \Z} \big] \\
        -\E \big[ \frac{\partial^2 f_1(\btheta, \Z)}{\partial \btheta^\top \partial \Z} - \frac{\partial^2 f_2(\btheta, \Z)}{\partial \btheta^\top \partial \Z} \big] & -\E \big[ \frac{\partial^2 f_1(\btheta, \Z)}{\partial \Z^\top \partial \Z} - \frac{\partial^2 f_2(\btheta, \Z)}{\partial \Z^\top \partial \Z} \big]
        \end{array}\right) 
\end{eqnarray*}
\end{comment}

where $\I_{11}$ and $\I_{22}$ are information matrix of $\btheta$ and $\Z$. We proceed to analyse the conditional and marginal distribution of $\widehat{\btheta}$. Theorem \ref{JointAsymDist} offers a strong foundation, from which corollaries follows, as Simulation \ref{Simulation:Beta-Bernoulli} illustrates.

\begin{corollary}{(Parameter Conditional Asymptotic Normality)}
    Suppose log ideal likelihood $h(\X, \Z|\btheta)$ and conditional probability function $f(\X |\Z, \btheta)$ satisfy Assumptions \ref{MarginalRegularity} and Assumption \ref{Negligibility} with some $N$. For invertible $\I_{11}$ and $\I_{22}$, the parameter marginal asymptotic distribution is
    \begin{equation}
        \sqrt{N}(\widehat{\btheta} - \btheta) | \Z \dconverge N\big(\textbf{0}, (\I_{11} - \I_{12}\I_{22}^{-1}\I_{12}^{\top})^{-1}\big),
        \label{MarginalTheta}
    \end{equation}
and conditional asymptotic distribution is 
    \begin{equation*}
        \sqrt{N}(\widehat{\btheta} - \btheta) | (\widehat{\Z}, \Z) \dconverge N\big(-\I_{11}^{-1} \I_{12}(\widehat{\Z}-\Z), \I_{11}^{-1}\big).
    \end{equation*}
    \label{ParameterConditionalAsymNorm}
\end{corollary}

See Supplement 3 for proof. \eqref{MarginalTheta} shows that, under conditional measure, MILE estimator $\widehat{\btheta}$ is asymptotically unbiased normal. Conditional on $\Z$ and its MILE estimator $\widehat{\Z}$, however, bias of $\I_{11}^{-1} \I_{12}(\widehat{\Z}-\Z)$ emerges. Thus, $\widehat{\btheta}$ should be viewed as estimators converging to pseudo true under the specific measure. Meanwhile, observe that $\I_{11}^{-1} \I_{12}(\widehat{\Z}-\Z) \pconverge 0$, asymptotically vanishing the bias. The corollary also apply to distribution of $\widehat{\Z}$ by swapping the roles of $\Z$ and $\bbeta$.

\subsection{Marginal Asymptotic Normality}
Up to this point, we developed the inference under conditional measures, leaving the problem under marginal measures. The empirical simulations shows that $\widehat{\btheta}$ may lose asymptotic normality marginally, but behave in terms of \eqref{TrueMarginalExpression} of Theorem \ref{TRUEMarginal}. Fortunately, normality is still verified under more conditions in \eqref{AsympTrueMarginalExpression}.

\begin{theorem}{(Marginal Asymptotic Distribution)}
    Suppose log ideal likelihood $h(\X, \Z|\btheta)$ and conditional probability function $f(\X |\Z, \btheta)$ satisfy Assumptions \ref{MarginalRegularity} with some $N$. There exists a positive definite matrix $\J_{\btheta, \Z}$, s.t. marginal asymptotic parameter distribution is
    \begin{equation}
        \widehat{\btheta} = \J^{-\frac{1}{2}}_{\btheta, \Z} \bmV + \btheta^c(\Z),
        \label{TrueMarginalExpression}
    \end{equation}
    where $\bmV$ is standard multivariate normal and $\btheta^c(\Z)$ is a random vector, for the estimator $\widehat{\Z} = \argmax_{\Z \in \mathcal{Z}} \ell(\Z, \btheta^{(0)} | \X)$, and some $\btheta^{(0)} \in \bTheta$.
    
$\widehat{\btheta}$ is asymptotic normal by Assumption \ref{Negligibility} with invertible $\I_{11}$ and $\I_{22}$, 
    \begin{equation}
        \sqrt{N}(\widehat{\btheta} - \btheta_c^{*}) \dconverge N(\mathbf{0}, \J^{-1}_{\btheta}),
        \label{AsympTrueMarginalExpression}
    \end{equation}
    for some $m_1 , m_2 \rightarrow \infty$, $N / m_1 \rightarrow 0$, $N / m_2 \rightarrow 0$, s.t. $\J^{-\frac{1}{2}}_{\btheta, \Z} = \O_p\Big(\frac{g_1(\btheta)}{N} + \frac{g_2(\btheta, \Z)}{m_1}\Big)$, $\btheta^c(\Z) = \btheta_c^{*} + \O_p\Big(\frac{g_3(\btheta, \Z)}{m_2}\Big)$, where $\btheta_c^{*}$ is constant, $g_1(\cdot)$, $g_2(\cdot, \cdot)$, $g_3(\cdot, \cdot)$ are bounded.
    \label{TRUEMarginal}
\end{theorem}

For brevity, denote $\sqrt{N} \J^{-\frac{1}{2}}_{\btheta, \Z} = \J^{-\frac{1}{2}}_{\btheta} + \O_p\Big(\frac{N g_2(\btheta, \Z)}{m_1}\Big)$.

Under conditions in Theorem \ref{TRUEMarginal}, which is mild, Simulation \ref{Simulation:Beta-Bernoulli} demonstrates marginal asymptotic normality. Additionally, MILE estimator is efficient and asymptotic unbiased with least variance who converges to Fisher's Information, $\lim_{N \rightarrow +\infty} \J_{\btheta} = \I(\btheta)$. Section \ref{sec:Simulation} presents simulation results, and Supplement 4 provides arguments of marginal normality, justifications of equivalent least variance, and a comparison between theoretical and empirical behaviours.

In practice, $\J_{\btheta}$ is computationally complicated. A feasible alternative is the Jackknife estimator $\widehat{\J_{\btheta}}$. As shown in \cite{Joe2005} and \cite{TF2024}, Jackknife delivers consistency in matrix form, $||\widehat{\J_{\btheta}} - \J_{\btheta}||_{\mathcal{F}} \pconverge 0$, under Frobenius norm.

\section{Simulation Studies}
\label{sec:Simulation}
We classify the common scenarios into two categories, differentiability with respect to $\Z$ and applicability of EM-type algorithm. Simulation studies are constructed to reflect the scenarios and benchmark performance against competitors.
\begin{enumerate}[label=(\alph*)~]
    \item Beta-Bernoulli mixture model, frequently employed in biostatistics data;
    \item Log-Cauchy mixture model, useful to capture heavy-tail pattern in financial data;
    \item Gaussian mixture model, widely-used to detect cluster;
    \item Bayesian Segmented Regression, applied in change point detection.
\end{enumerate}

In general, we notate $N$ as the number of units/individuals  characterized by latent effects $z_i$. Similarly for individual $i$, $x_{ij}$ is the $i$-th observation, where $j = 1, 2, \cdots, M$. These notations are uniformly applied to examples in this Section.

\subsection{Beta-Bernoulli Mixture Model}
\label{Simulation:Beta-Bernoulli}
Simulation 1 considers the Beta-Bernoulli Mixture Model. We assume that the observations $x_{ij}$ follows Bernoulli distribution conditional on latent random effects, $x_{ij}|z_{i} \sim Bern(z_i)$, where $z_i$ follows Beta distribution $z_i \sim Beta(\alpha, \beta).$

We set $\alpha = \beta = \theta$, and consider $\theta = 5$ and $\theta = 10$ in simulations. Individual number is taken as $N \in \{10, 20, 50, 100, 200, 500, 1000\}$, with $M \in \{1000, 10000\}$ observations per individual. Implementations of the EM and MILE estimators are presented in Supplement 5. Simulation outcomes are reported in Table \ref{BetaBerTable}. The latent estimator $\widehat{z}$ performs well, with small bias and stable standard deviation. 

For large $M$ and $N$, $\widehat{\theta}_{MILE}$, the MILE estimator $\widehat{\theta}_{MILE}$ shows excellent agreement with the EM estimator $\widehat{\theta}_{EM}$, with almost same standard deviations. Since EM estimator is consistent to MLE and achieves the asymptotic efficiency, the results indicate that MILE estimator empirically shares the same asymptotic variance as the MLE. Supplement 4 contains more verifications.

\subsection{Log-Cauchy Mixture Model}
\label{Simulation:logCauchy}

In Simulation \ref{Simulation:logCauchy}, we consider $N$ investors, each holding $M$ pension fund shares. $x_{ij}$ is the discounted beneficiable value of $i$-th individual's $j$-th share. Conditional on the investor lifetimes, i.e. latent effect $z_i$, we model $x_{ij}|z_i \sim N(e^{-rz_i}, \sigma_1^2)$, and assume $z_i$ follows log-Cauchy distribution $z_i \sim logCauchy(\mu, \sigma_2^2)$. Set $\mu = 2$, $\sigma_1 = \sigma_2 = 1$, $r = 0.05$, with $M$ and $N$ matching Simulation \ref{Simulation:Beta-Bernoulli}. Results are presented in Table \ref{logCauchyTable}.

In Simulation \ref{Simulation:logCauchy}, the EM algorithm cannot be applied, so the moment estimators (MoM) are used as competitors instead. See Supplement 5 for reasoning. MILE outperforms MoM, exhibiting reduced variability (``Sd" rows) under moderate and small sample sizes. Moreover, the latent variables are estimated directly, accurately and with smaller deviation.

\subsection{Gaussian Mixture Model (GMM)}
\label{Simulation:GMM}

GMM is a well-defined model with well-known structure. Set sample size to be $N$ with $M = 1$ and category number $K = 3$. SCP is a applied to GMM scenario. Results are summarised in Table \ref{GMMTable}. 

\begin{figure}[hbtp]
    \centering
    \includegraphics[width=0.78\linewidth]{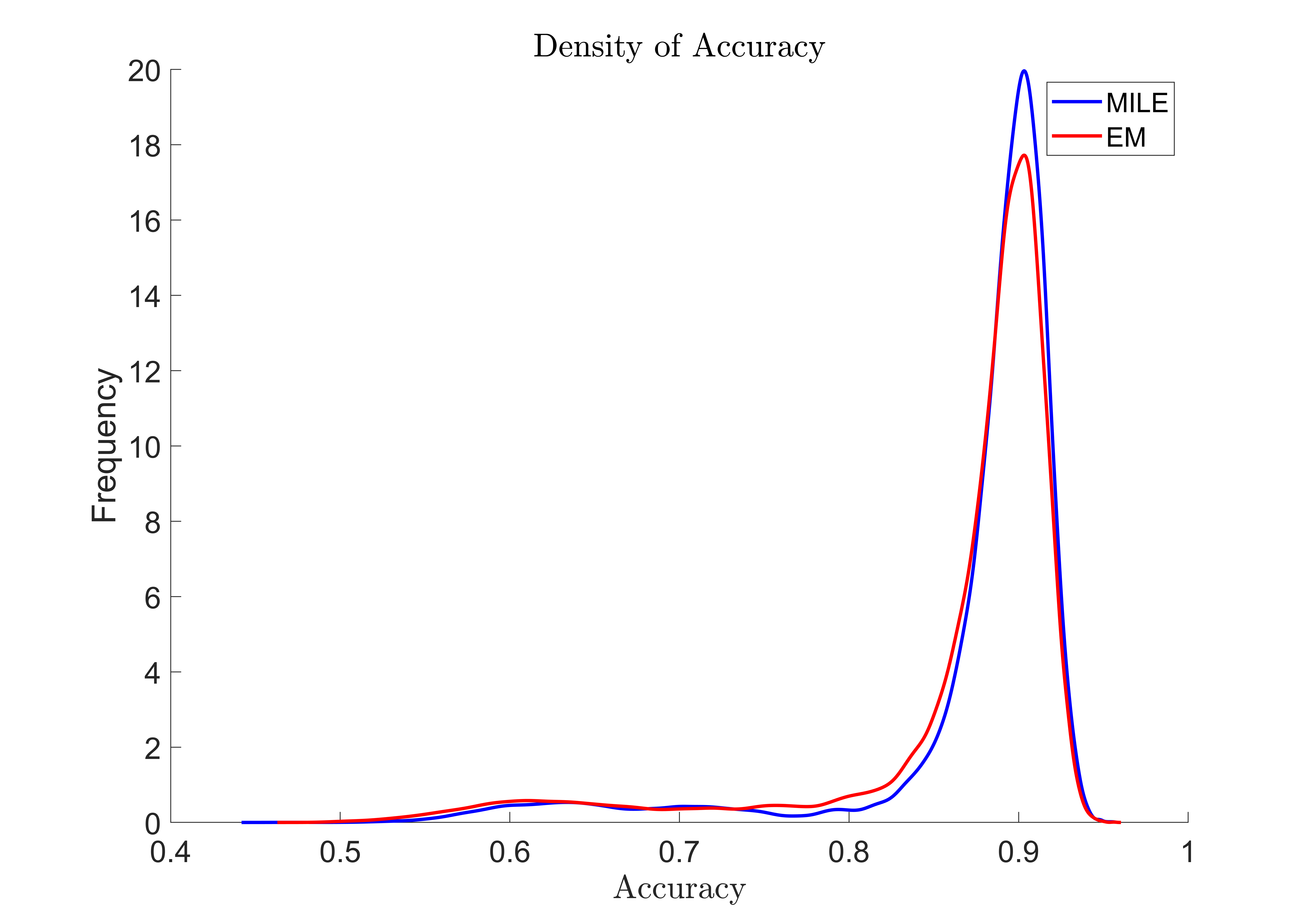}
    \caption{Density Plot, Accuracy of MILE and the EM prediction}
    \label{GMMAccuDensity}
\end{figure}

While MILE and EM estimators exhibit comparable bias and standard deviations, MILE dominantly outperforms the EM algorithm in accuracy, with higher mean but less standard deviation. Figure \ref{GMMAccuDensity} highlights the pattern that MILE produces more high-accuracy, but fewer low-accuracy replicates than the EM algorithm. GMM results are used as SCP initials, but MILE still converges rapidly.

\subsection{Bayesian Segmented Regression}
\label{Simulation:BayesianSegmentRegression}

\begin{figure}[hbtp]
    \centering
    \includegraphics[width=0.65\linewidth]{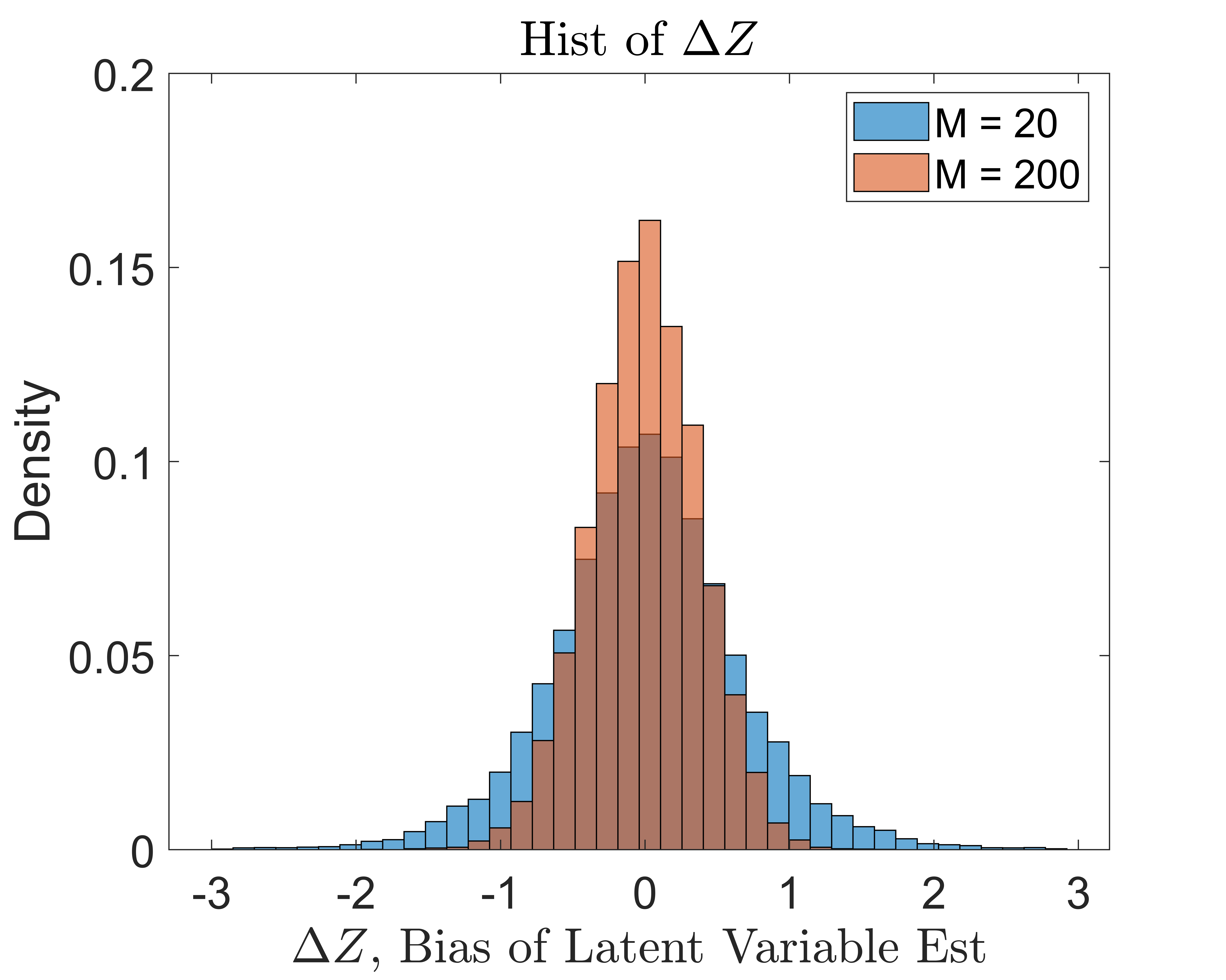}
    \caption{Histogram of the Bias of Change Point Estimates, $N$ = 10}
    \label{BSRZhist}
\end{figure}

In Simulation \ref{Simulation:BayesianSegmentRegression}, we consider $N$ independent time series, each observed at the same set of sampled timestamp, $t \in \{t_1, t_2, \cdots, t_M\}$, where $0 < t_j < T$. Each time series has single but different change point, occurring at $z_i$ for series $i$, and $x_{i,t}$ is observations at time $t$. Change points have a scaled Beta prior, $z_i/T \sim Beta(\alpha, \beta)$. Hierarchically, observations follow independent heterogeneous Poisson distribution,
\begin{equation}
    x_{i,t} \sim
    \begin{cases}
      Pois(\lambda^1_{i, t}) & t < z_i, \\
      Pois(\lambda^2_{i, t}) & t \geq z_i,
    \end{cases}
\end{equation}
where $\lambda^1_{i, t}=e^{\beta_1\left(t-z_i\right)+a}$ and $\lambda^2_{i, t}=e^{\beta_2\left(t-z_i\right)+a}$. Results are shown in Table \ref{BSRTable}.

\begin{figure}[hbtp]
    \centering
    \begin{minipage}{0.45\textwidth}
        \centering
        \includegraphics[width=0.85\textwidth]{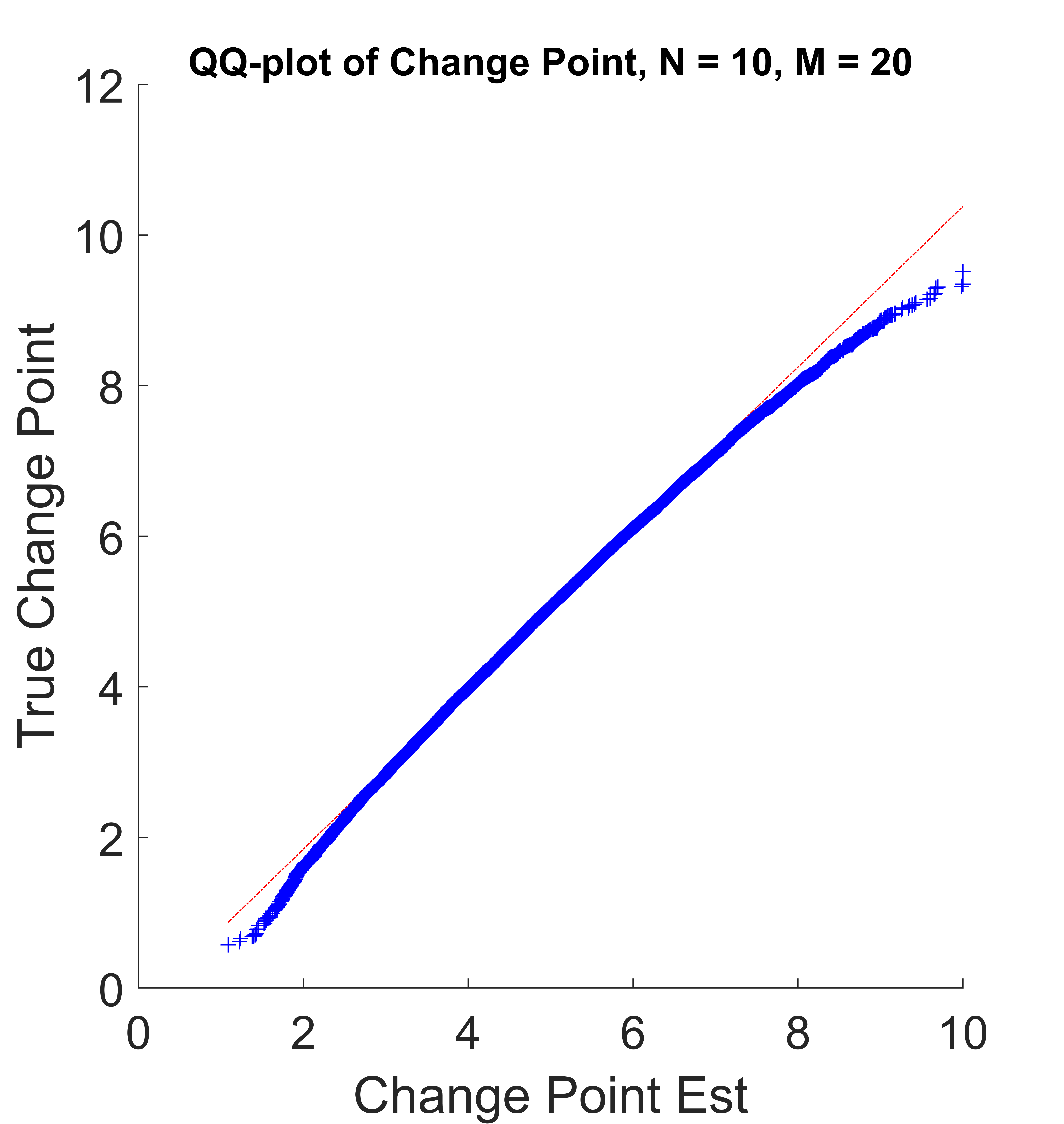}
        \caption{$M$ = 20, $N$ = 10}
        \label{BSRQQplot20}
    \end{minipage}
    \begin{minipage}{0.45\textwidth}
        \centering
        \includegraphics[width=0.85\textwidth]{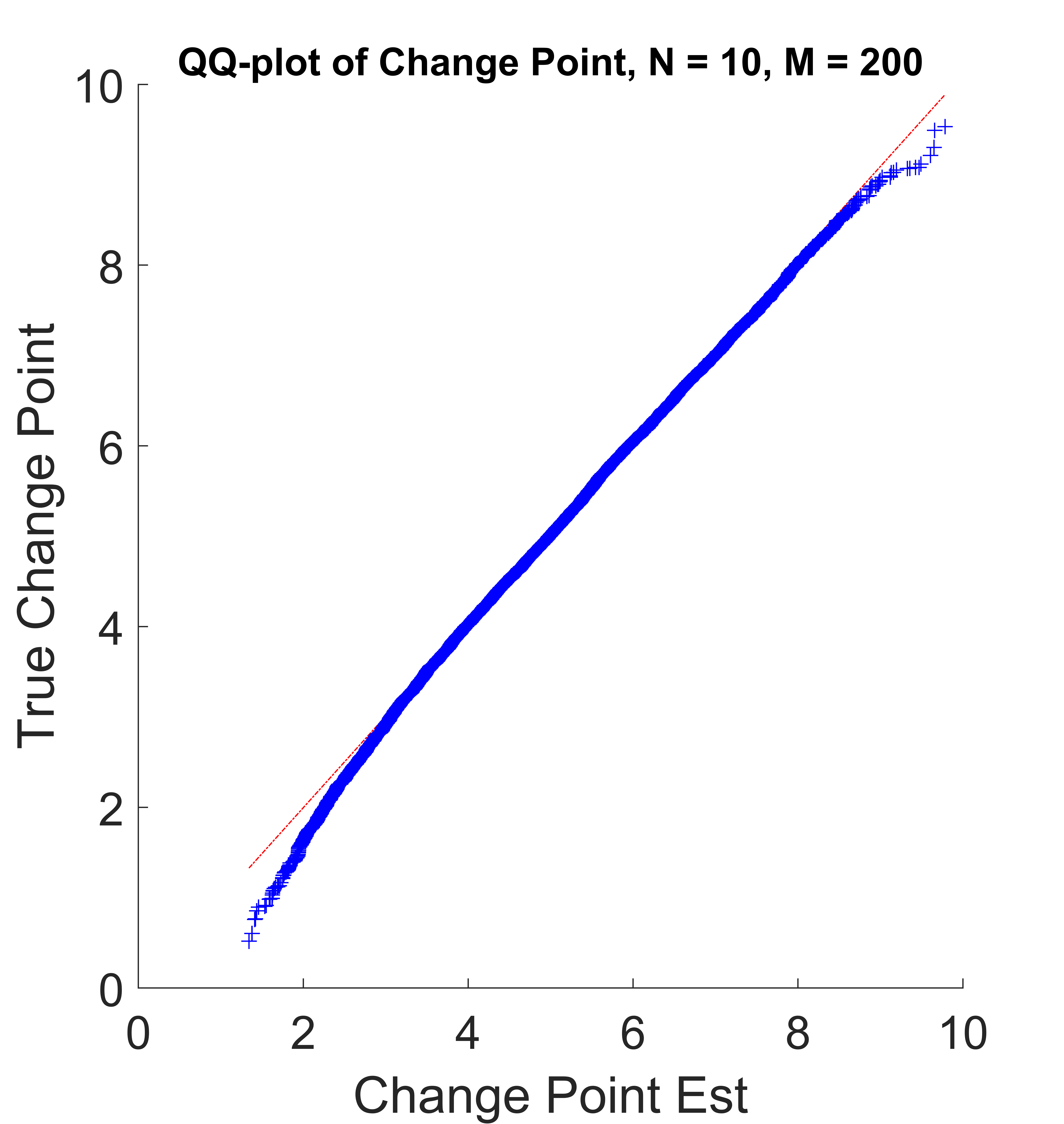} 
        \caption{$M$ = 200, $N$ = 10}
        \label{BSRQQplot200}
    \end{minipage}
\end{figure}

Change points of BSR are the primary focus, and thus we evaluate latent variable estimates. Figure \ref{BSRZhist} shows that the accuracy of estimators $\widehat{\Z}$ improves under larger $M$, consistent with the theoretical results in Section \ref{sec:Inference}. Empirically, the standard deviation of $\widehat{\Z}$ in Table \ref{BSRTable} decreases under large sample size. QQ plots in Figures \ref{BSRQQplot20} \ref{BSRQQplot200} indicate the empirical distribution of $\widehat{\Z}$ closely approximate normality.

\section{Illustrative Data Example} \label{sec:data}

We analyse a environmental dataset, to illstrate the proposed framework, examining the relationship between Land Surface Temperature (LST) Data and ecoregions in United States. LST, measured separately for daytime (DLST) and night (NLST, reflects thermal energy flow among land surface, atmosphere, and biosphere, thereby indicating local environmental conditions \citep{Ma2023JABES}. Although influenced by multiple factors, LST is strongly associated with air temperature, which is increasing due to global warming \citep{NOAA2021, IPCC2021}.

Ecoregion is a expertise-selected class of climate categories that represents ecosystem areas with similarity, including the type and quantity of environmental resources. The ecoregion framework, originally developed by \cite{Omernik1987} and further refined through collaborative mapping efforts with Environmental Protection Agency (EPA) regional offices and other agencies, provides a spatial foundation for ecosystem research and assessment. Ecoregions identify similar areas that play a vital role in guiding ecosystem management strategies, which are often responsible for managing different resources within the same geographic regions \citep{Omernik2014, McMahon2001}. Ecoregions are hierarchical into four levels: low levels capture broad patterns and high levels distinguish specific units.

In practice, We use LST data with additional covariates to predict ecoregions across the United States. The covariates include spatial coordinate (latitude, longitude) and soil, vegetation and hydrology characteristics. Let there be $n$ locations and $K$ ecoregions. If location $i$ belongs to ecoregion $k$, the latent cluster label , where $\mathbbm{P}(Z_i = k) = p_k$ and $\sum_{k = 1}^K p_k = 1$. Conditional on $Z_i = k$, observations $\Y_i$ follows $\Y_i|(Z_i = k) \sim N(\bmu_{i, k}, \bSigma_k)$ independently, where $\bmu_{i, k} = (\X_i \bbeta_{k, 1}, \X_i \bbeta_{k, 2})$ and $\X_{i}$ denotes covariates of location $i$.

\begin{figure}[hbtp]

\caption{\scriptsize Evolution of cluster assignments for U.S. LST Data under SCP. Each panel displays current grouping as latent variables. Panels are arranged from left to right, top to bottom, with 2000 iterations between adjacent panels.}
    \centering
    % First row
    \begin{minipage}{0.32\textwidth}
        \centering
        \includegraphics[width=\linewidth]{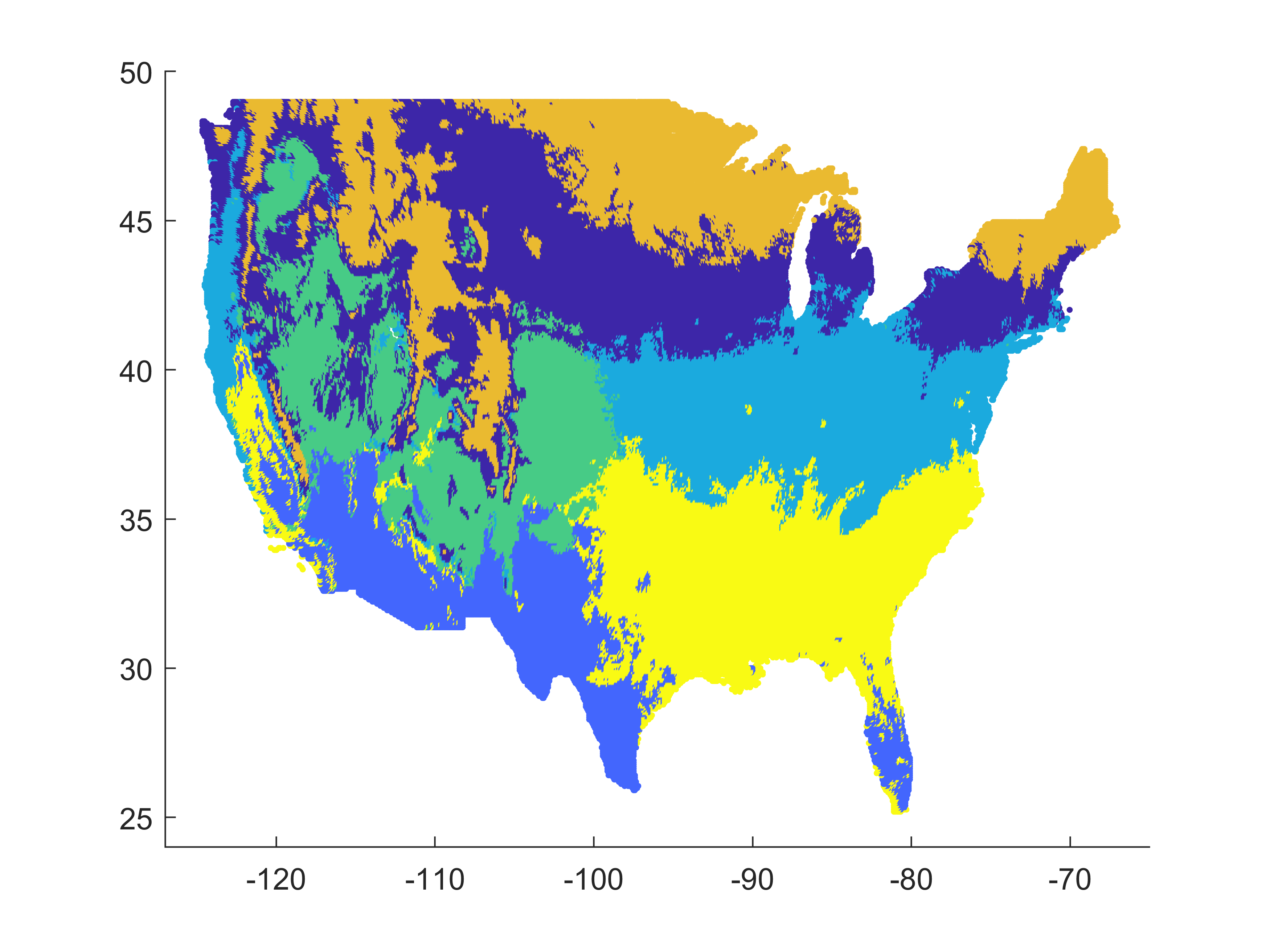}
    \end{minipage}
    \hfill
    \begin{minipage}{0.32\textwidth}
        \centering
        \includegraphics[width=\linewidth]{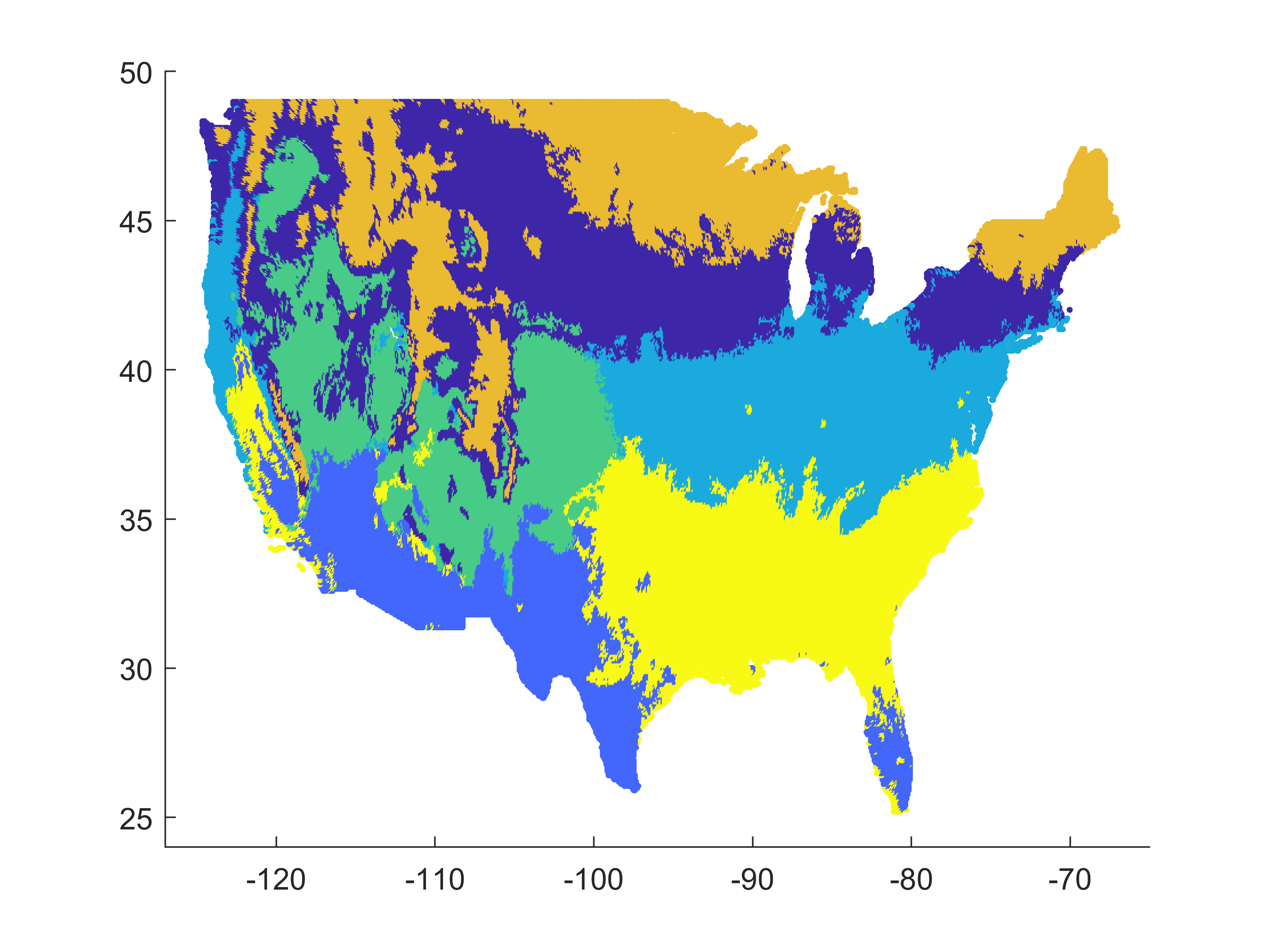}
    \end{minipage}
    \hfill
    \begin{minipage}{0.32\textwidth}
        \centering
        \includegraphics[width=\linewidth]{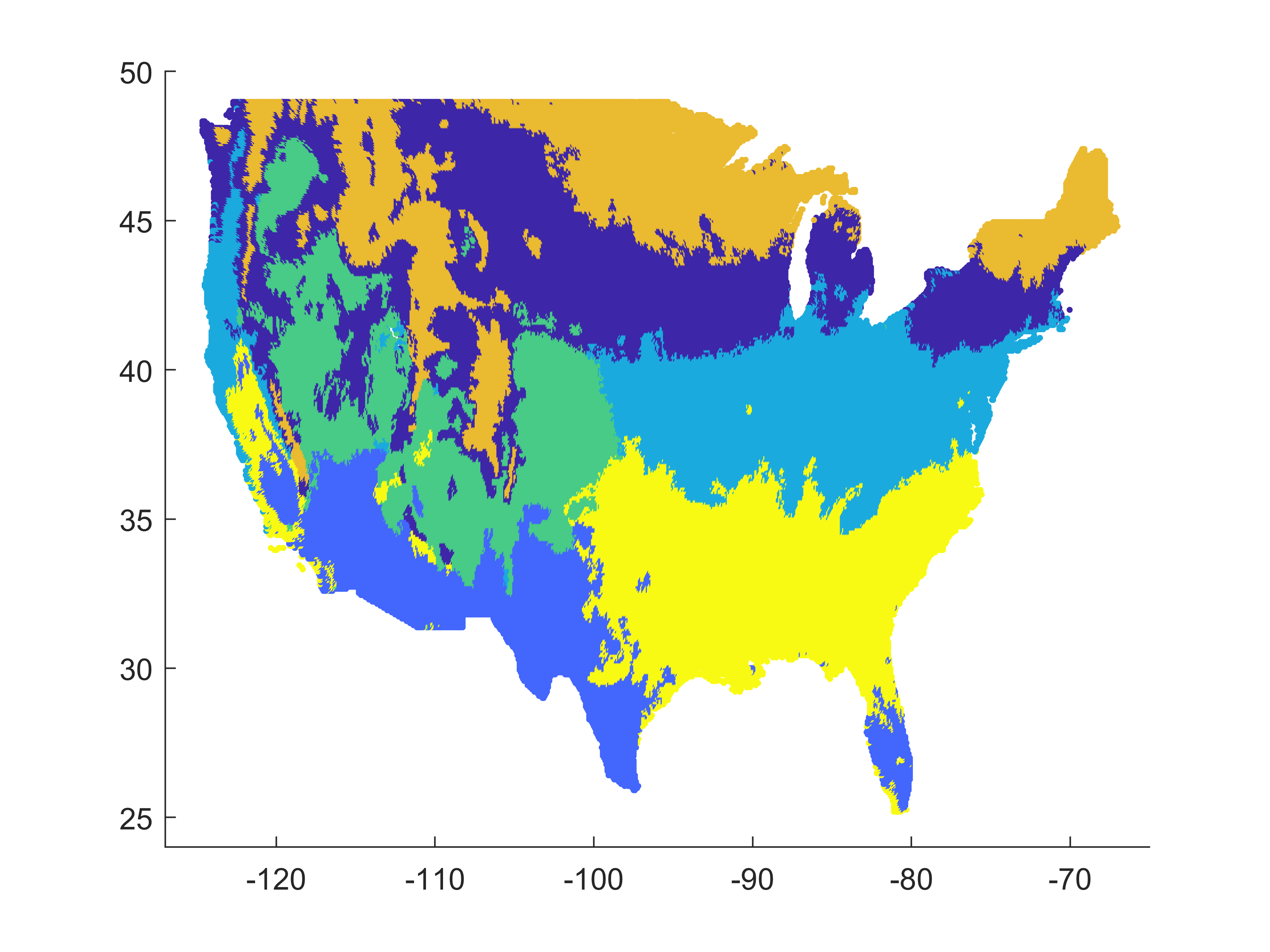}
    \end{minipage}
    % Second row
    \begin{minipage}{0.32\textwidth}
        \centering
        \includegraphics[width=\linewidth]{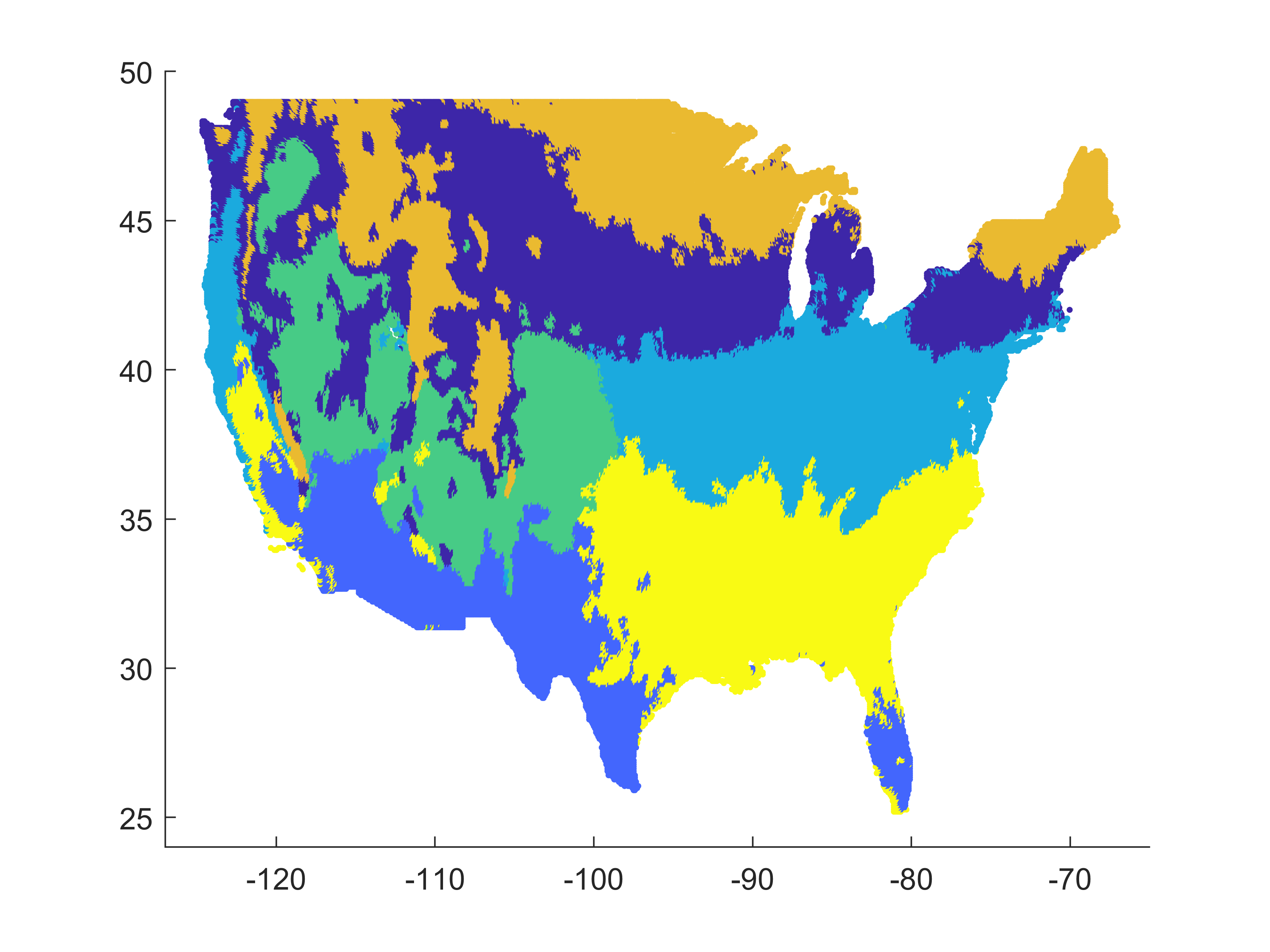}
    \end{minipage}
    \hfill
    \begin{minipage}{0.32\textwidth}
        \centering
        \includegraphics[width=\linewidth]{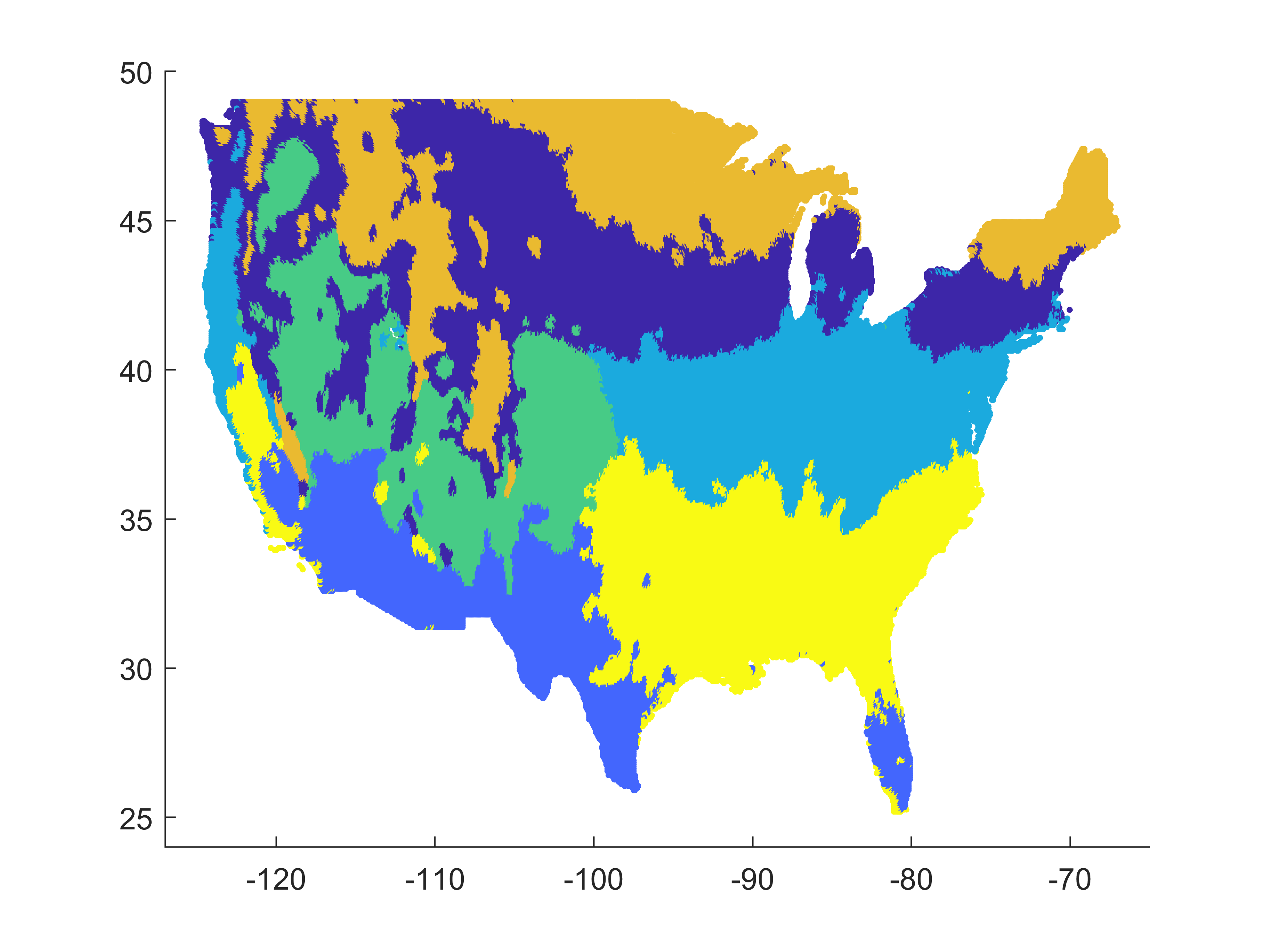}
    \end{minipage}
    \hfill
    \begin{minipage}{0.32\textwidth}
        \centering
        \includegraphics[width=\linewidth]{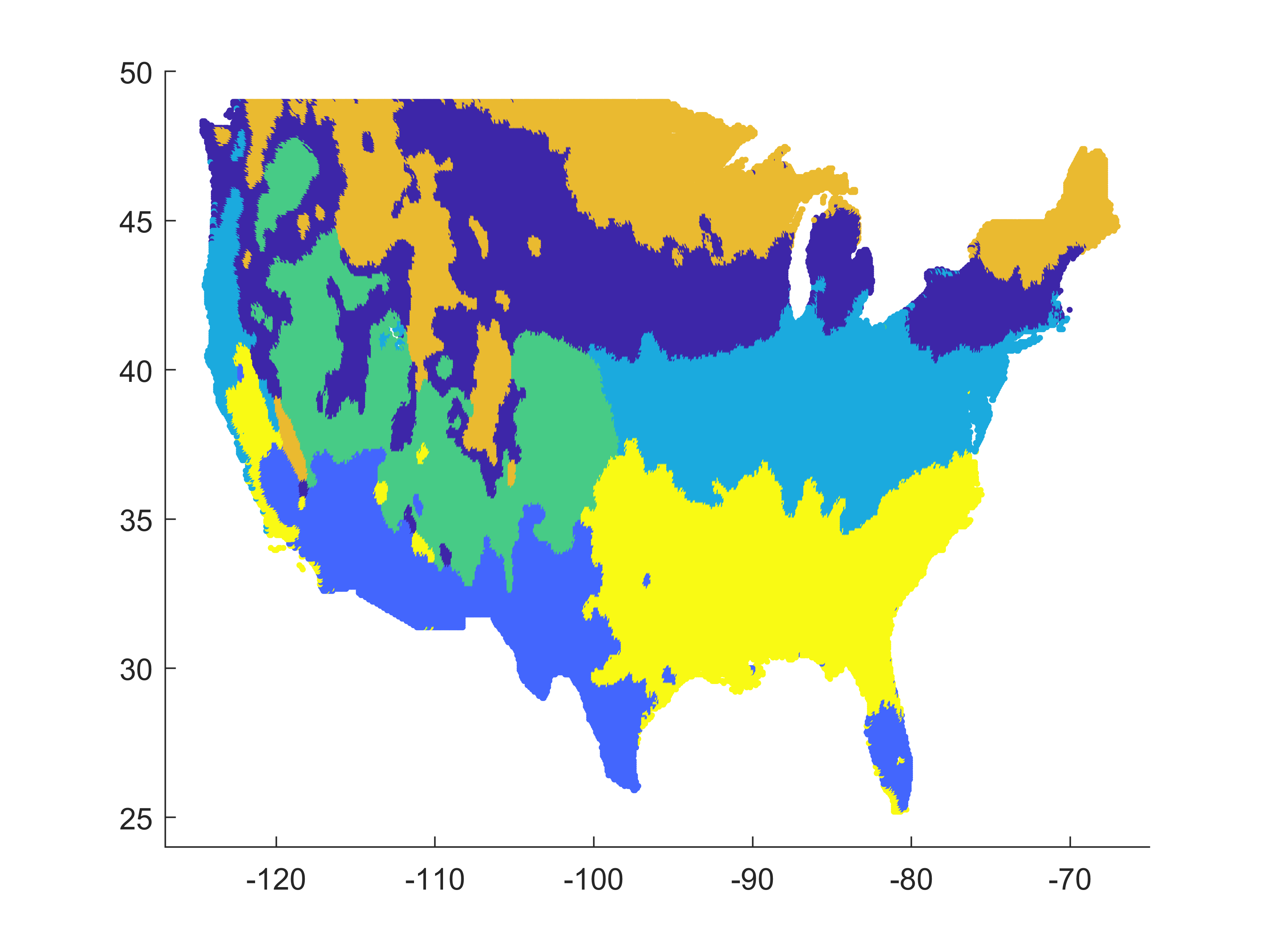}
    \end{minipage}
    % Third row
    \begin{minipage}{0.32\textwidth}
        \centering
        \includegraphics[width=\linewidth]{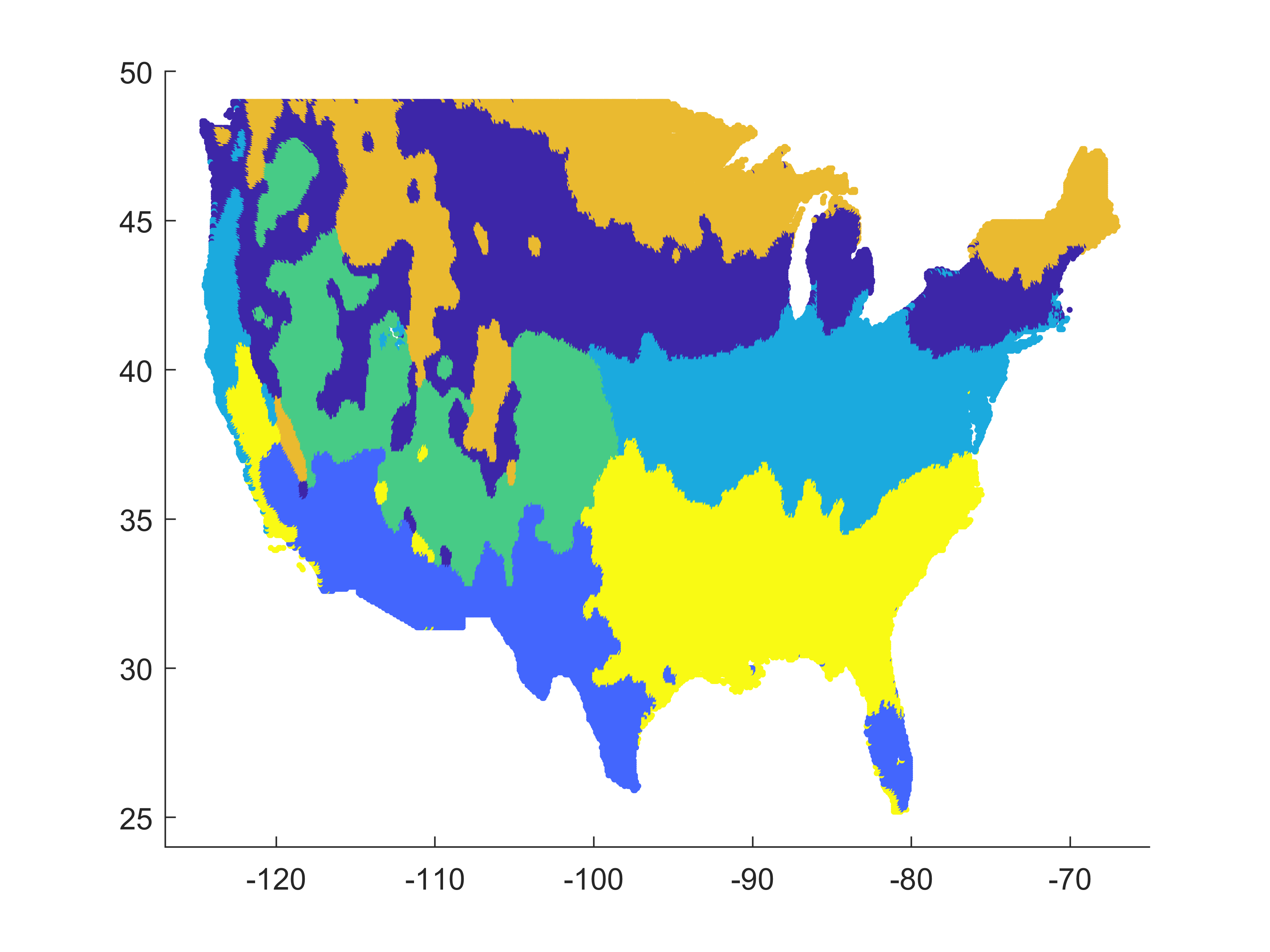}
    \end{minipage}
    \hfill
    \begin{minipage}{0.32\textwidth}
        \centering
        \includegraphics[width=\linewidth]{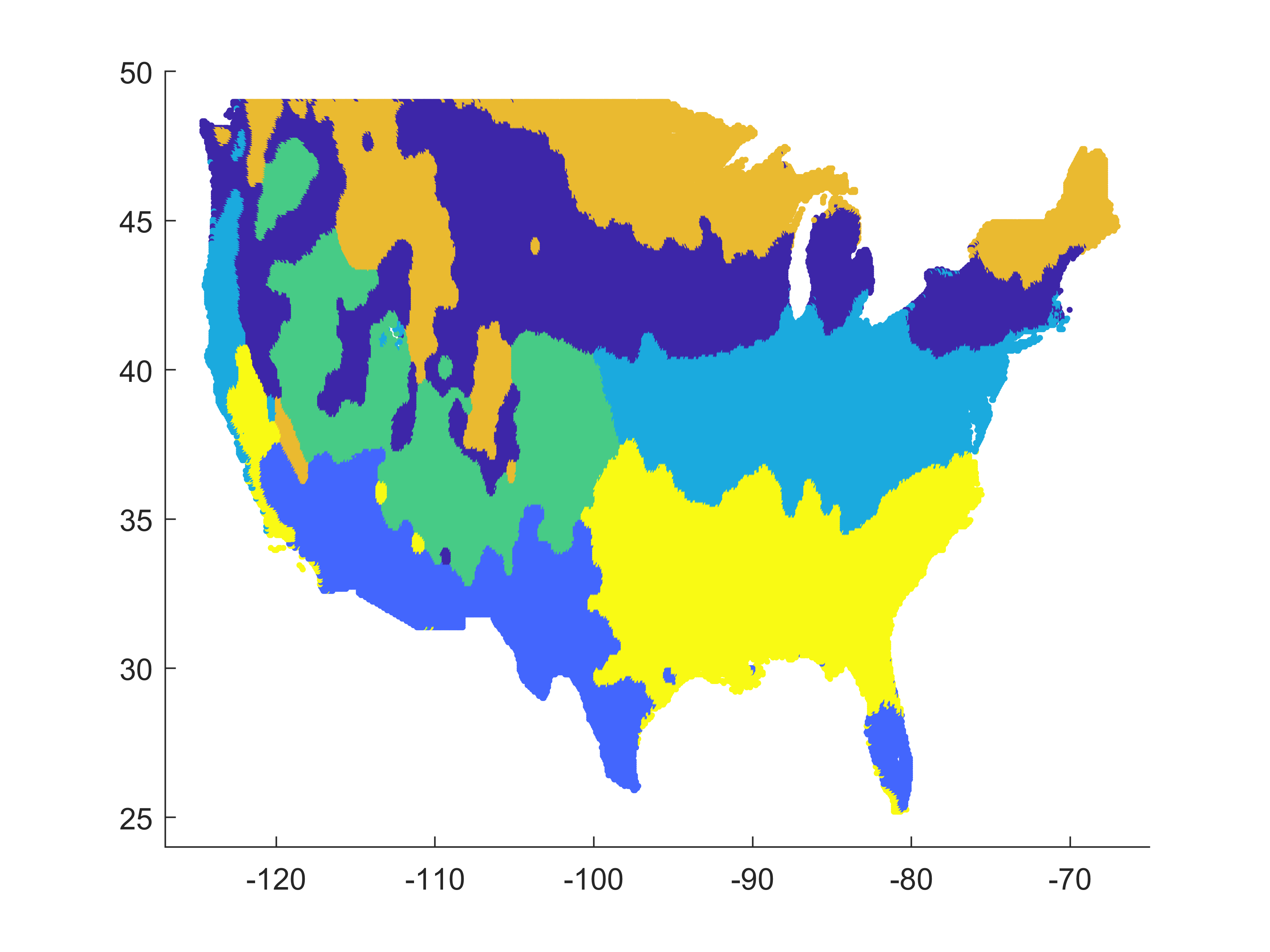}
    \end{minipage}
    \hfill
    \begin{minipage}{0.32\textwidth}
        \centering
        \includegraphics[width=\linewidth]{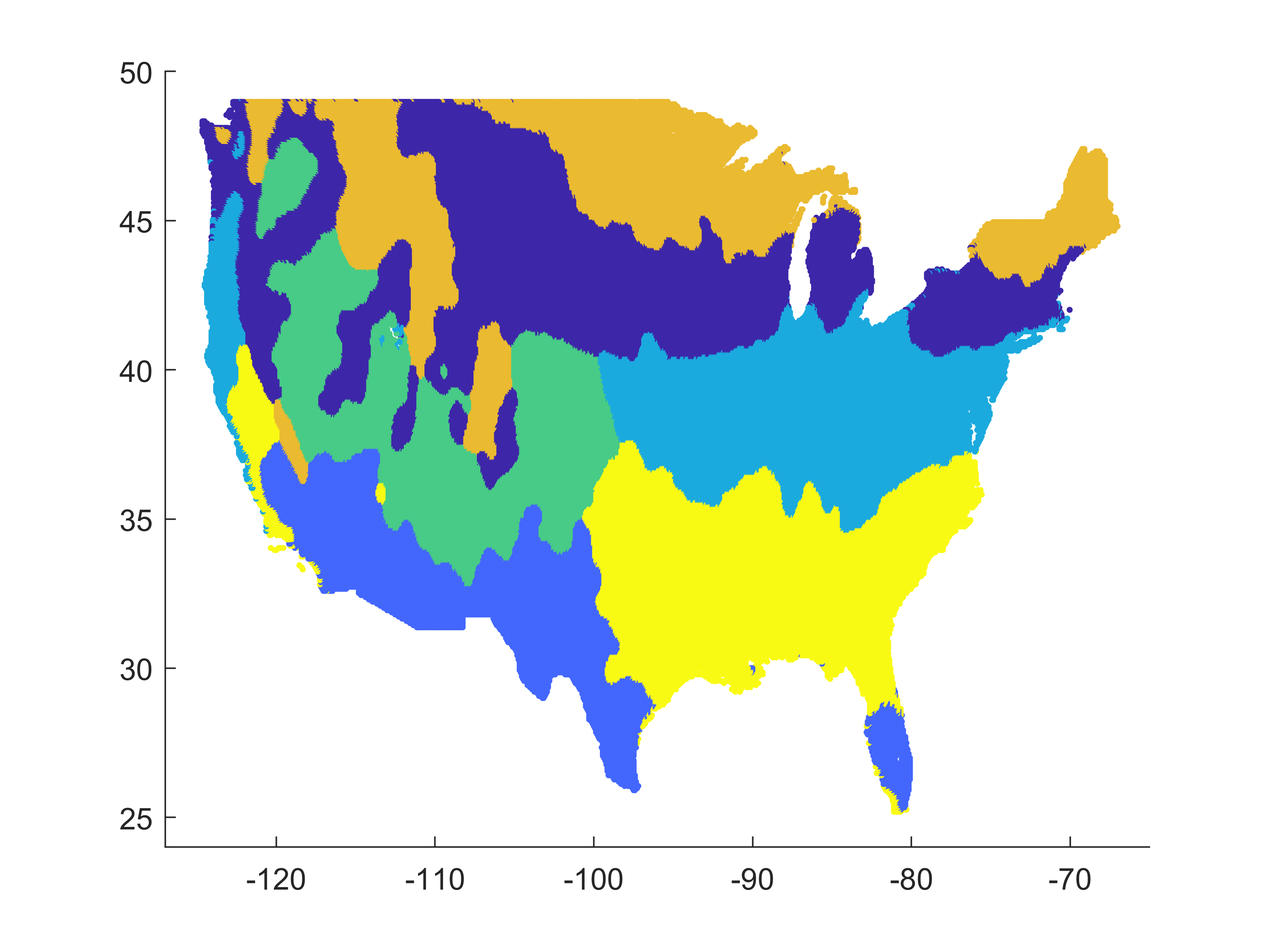}
    \end{minipage}

\label{LSTMAP}
\end{figure}

Since Level 1 ecoregion of Unites State contains 12 categories, many of which correspond to very small land area, we set $K = 6$ to reflect the major zones. MILE is estimated via SCP, initialized by multi-dimensional k-means clustering. Figure \ref{LSTMAP} shows the improvement trajectory and cluster assignment overlaid onto the map. 

\section{Discussion and Future Work}
\label{sec:Discussion}
MILE provides a unified framework for latent variable models by maximizing the ideal likelihood, thereby enabling simultaneous estimates of parameters and latent variables. The framework is broadly applicable and robust, particularly in extreme scenarios where traditional methods fails. Similar to the EM algorithm, MILE achieves asymptotic efficiency with least variances, but outperforms other competitors such as MCMC and MoM. Simulation studies and theoretical results, including asymptotic properties, under mild conditions support its reliability, shaping MILE a comprehensive alternative to existing methods.

\subsection{Approximation Compatibility}

Flexibility of MILE to target functions extends itself to likelihood or distribution approximation. Suitable assumptions could be imposed to ensure empirical validity, permitting implementation of MILE under approximations in various settings.

In distribution approximation, well-established methods provide inference with convincing results. For instance under Bayesian structure, Variational Bayes \citep{Jordan1999} seeks approximation $\widehat{Q}(\Z | \btheta)$ in a distribution family $\Q$, minimising the K-L divergence, to the posterior $f(\Z | \btheta)$. See \cite{VBReview} for reviews. In MILE framework, the posterior distribution is replaced by ideal likelihood as

$$
\widehat{Q}(\Z | \btheta) =\argmax_{\btheta \in \Theta} \argmin_{Q(\Z | \btheta) \in \Q}  \int_{\mathcal{Z}} \log \frac{Q(\Z | \btheta)}{f(\X, \Z | \btheta)} dQ(\Z | \btheta).
$$

Integrated Nested Laplace Approximation (INLA) provides another perspective on distribution approximation. Originally developed for the class of latent Gaussian models (LGM), it employs Laplace's Method as an accurate and fast alternative to MCMC \citep{INLA2009}. It is natural to replace the target distribution by ideal likelihood from the view of MILE. Extension beyond the LGM settings have been studied \cite{RUE2014, INLANonLGM2024}, with comprehensive inference by \cite{INLAbook2020} and Bayesian asymptotics by \cite{Miller2021}. 

Approximate likelihood provides practical solutions when the full likelihood is intractable. Notable approaches, including Composite Likelihood \citep{Lindsay1988, Varin2011} and Vecchia approximation \citep{Vecchia1988,VecchiaReview2021}, achieve scalability in spatial statistics studies. Replacing the full likelihood by ideal likelihood, MILE inherits advantages under mild conditions.

\begin{comment}
Approximate methods usually readily reduce to optimisation problems. Moreover, strict global convexity is usually assumed in literature for statistical inference. Strict global convexity implies unique extremum on the support, guaranteeing numerical fast convergence of many optimisation routines. For comprehensive illustrations, please refer to \cite{cvx2004} for convex optimisation, and \cite{wright2006} for numerical solutions and its analysis.    
\end{comment}

While differentiability is assumed in optimisation problems, no-gradient problems remain significant. These targets necessitate powerful derivative-free methods, aligned with MILE applications with discrete latent variables. Representative approaches include toolbox of \cite{Liu2022} and discrete optimisation methods by \cite{Choirat2012}, \cite{Kozek1998} and \cite{Ma2023}.

Overall, approximations and other numerical methods could be integrated with MILE to overcome numerical challenges. Incorporating these techniques, MILE retains flexibility and yield reliable estimations in a wide range of applications.

\subsection{Algorithm \& Inference Improvement}
GA for latent values could be improved regarding convergence and speed. \cite{Richter2017} and \cite{Tian2009} proposed their work that help GA improvements. We left the promising directions in future.

A further inquiry concerns inference under latent variable estimates. Partitioning parameters into $\btheta = (\balpha, \bbeta)$, as presented in many models, yields factorization $L(\btheta, \Z; \X) = f(\X|\Z, \bbeta) f(\Z|\balpha)$. Under the assumption, estimator $\widehat{\balpha}$ are determined by $\widehat{\Z}$, and its consistency relies on the accuracy of $\widehat{\Z}$. Given rapid convergence of $\widehat{\Z}$ to $\Z$, $\widehat{\balpha}$ is also consistent. Whereas, slow convergence may lead to $\widehat{\balpha}$ converging to a pseudo true; see Simulation \ref{Simulation:BayesianSegmentRegression}. Treating $\widehat{\Z}$ as observed $\Z$ with error links the problem to $\widehat{\balpha}$ under measurement error, where bias depends on specific structure of measurement error. Future research should clarify when parameters converge to true value under what conditions, regardless of the convergence rate.

\section{Tables}

% Table generated by Excel2LaTeX from sheet 'ComparisonTable'
\begin{table}[hbtp]
  \centering
  \small
  \caption{Comparison Table}
    \begin{tabular}{llllllll}
    \toprule
    \toprule
    $\pi_{\btheta}$ & $\partial_{\Z}$& $\Z_{\bpi^*}$& $\mathbbm{E}_{\Z_{\bpi^*}}$ & (MC)EM    & MCMC  & MILE & Speed \\
    
    \ding{51}     & $\forall$     & \ding{51}     & \ding{51}     & \ding{51}    & \ding{51}     & \ding{51}     & ? \\
    \ding{53}     & $\forall$     & $\forall$     & $\forall$     & ?    & \ding{53}     & \ding{51}     & ? \\
    ?    & \ding{51}     & \ding{51}     & \ding{51}     & \ding{51}     & \ding{53}     & \ding{51}     & $\gtrdot$ \\
    ?    & \ding{51}     & \ding{53}     & \ding{51}     & \ding{51}     & \ding{53}     & \ding{51}     & $\gtrdot$ \\
    ?    & \ding{51}     & \ding{51}     & \ding{53}     & \ding{53}     & \ding{53}     & \ding{51}     & $\gtrdot$ \\
    ?    & \ding{51}     & \ding{53}     & \ding{53}     & \ding{53}     & \ding{53}     & \ding{51}     & $\gtrdot$ \\
    ?    & \ding{53}     & \ding{51}     & \ding{51}     & \ding{51}     & \ding{53}     & \ding{51}     & $\lessdot$ \\
    ?    & \ding{53}     & \ding{53}     & \ding{51}     & \ding{53}     & \ding{53}     & \ding{51}     & $\lessdot$ \\
    ?    & \ding{53}     & \ding{51}     & \ding{53}     & \ding{53}     & \ding{53}     & \ding{51}     & $\lessdot$ \\
    ?    & \ding{53}     & \ding{53}     & \ding{53}     & \ding{53}     & \ding{53}     & \ding{51}     & $\lessdot$ \\
    \toprule
    \toprule
    \end{tabular}%
  \label{ComparisonTable}%
\end{table}%

%%%%%%%%%%%%%%
% Table generated by Excel2LaTeX from sheet 'BetaBerToy'

\begin{table}[hbtp]
  \centering
  \scriptsize
  \caption{\scriptsize Results for Simulation 1 (Beta-Bernoulli Mixture Model). $M$ is numbers of observations per individual and $N$ is number of individuals. ``Est" and ``Sd" represent mean and standard deviation of estimator $\widehat{\btheta}$ over 5000 Monte Carlo experiments. ``time" records average time cost (in second) over 100 Monte Carlo experiments. Latent variable estimation is evaluated by the average bias $\Delta Z = mean(\widehat{Z} - Z)$ and its standard deviation ``Sd, $\Delta Z$".}
  \begin{tabular}{clrrrrrrrrrr}
    \toprule
    \toprule
       &    & \multicolumn{5}{c}{$M$ = 1000, True = 5} & \multicolumn{5}{c}{$M$ = 1000, True = 10} \\
    \multicolumn{2}{c}{$N$} & 10 & 20 & 50 & 100 & 1000 & 10 & 20 & 50 & 100 & 1000 \\
    \multirow{3}[2]{*}{EM} & Est & 6.71 & 5.36 & 5.03 & 5.07 & 4.97 & 12.18 & 10.90 & 10.30 & 10.22 & 10.02 \\
       & Sd & 4.16 & 1.58 & 0.91 & 0.69 & 0.21 & 6.27 & 3.82 & 2.08 & 1.48 & 0.45 \\
       & time & 1.60 & 1.59 & 2.12 & 2.46 & 9.82 & 1.98 & 2.98 & 4.00 & 4.43 & 18.62 \\
    \multirow{5}[2]{*}{MILE} & Est & 6.42 & 5.52 & 5.21 & 5.15 & 5.04 & 12.42 & 11.25 & 10.58 & 10.42 & 10.18 \\
       & Sd & 4.50 & 1.67 & 0.99 & 0.73 & 0.21 & 7.02 & 3.94 & 2.26 & 1.54 & 0.46 \\
       & $\Delta$Z & -0.00 & 0.00 & 0.00 & 0.00 & 0.00 & -0.00 & 0.00 & 0.00 & 0.00 & 0.00 \\
       & Sd, $\Delta$Z & 0.01 & 0.02 & 0.02 & 0.01 & 0.02 & 0.02 & 0.02 & 0.02 & 0.02 & 0.02 \\
       & time & 1.90 & 1.98 & 2.13 & 2.25 & 5.27 & 2.08 & 2.07 & 2.17 & 2.26 & 4.34 \\
       &    & \multicolumn{5}{c}{$M$ = 10000, True = 5} & \multicolumn{5}{c}{$M$ = 10000, True = 10} \\
    \multirow{3}[2]{*}{EM} & Est & 5.77 & 5.41 & 5.12 & 5.02 & 5.01 & 12.00 & 11.05 & 10.48 & 10.15 & 10.06 \\
       & Sd & 2.70 & 1.87 & 0.91 & 0.67 & 0.19 & 5.55 & 3.75 & 2.01 & 1.32 & 0.43 \\
       & time & 2.47 & 4.39 & 6.41 & 10.17 & 152.12 & 3.30 & 5.65 & 7.92 & 11.73 & 277.59 \\
    \multirow{5}[2]{*}{MILE} & Est & 6.23 & 5.56 & 5.23 & 5.11 & 5.01 & 12.44 & 10.99 & 10.40 & 10.19 & 10.02 \\
       & Sd & 3.43 & 1.88 & 1.06 & 0.73 & 0.21 & 6.42 & 3.66 & 2.14 & 1.45 & 0.43 \\
       & $\Delta$Z & 0.00 & 0.00 & 0.00 & 0.00 & 0.00 & 0.00 & 0.00 & 0.00 & 0.00 & 0.00 \\
       & Sd, $\Delta$Z & 0.00 & 0.00 & 0.00 & 0.00 & 0.00 & 0.00 & 0.00 & 0.00 & 0.00 & 0.00 \\
       & time & 2.25 & 2.90 & 3.51 & 4.86 & 24.47 & 1.90 & 2.42 & 3.22 & 4.32 & 22.94 \\
    \bottomrule
    \bottomrule
  \end{tabular}
  \label{BetaBerTable}
\end{table}

\begin{table}[hbtp]
  \centering
  \scriptsize
  \caption{\scriptsize Results for Simulation 2 (log-Cauchy Mixture Model). $M$ is numbers of observations per individual and $N$ is number of individuals. ``Est" and ``Sd" represent mean and standard deviation of estimator $\widehat{\btheta}$ over 5000 Monte Carlo experiments. ``time" records average time cost (in second) over 100 Monte Carlo experiments. Latent variable estimation is evaluated by bias median $\Delta Z_{(m)}$" and its standard deviation of bootstrap sample median ``Sd, $\Delta Z_{(m)}$".}
  \begin{tabular}{clrrrrrrrrrrrr}
    \toprule
    \toprule
    & & \multicolumn{5}{c}{$M=1000$, True = 2} & & \multicolumn{5}{c}{$M=10000$, True = 2} \\
    \multicolumn{2}{c}{$N$} & 10 & 20 & 50 & 100 & 1000 & & & 10 & 20 & 50 & 100 & 1000 \\

    \multirow{3}{*}{MoM} 
       & Est & 1.87 & 1.91 & 2.02 & 1.97 & 1.99 & & & 1.93 & 1.94 & 1.98 & 1.97 & 2.00 \\
       & Sd & 0.68 & 0.38 & 0.26 & 0.21 & 0.06 & & & 0.73 & 0.45 & 0.26 & 0.20 & 0.06 \\
       & time & 0.02 & 0.02 & 0.02 & 0.03 & 2.20 & & & 0.02 & 0.03 & 0.07 & 0.23 & 18.33 \\
    \multirow{5}{*}{MILE} 
       & Est & 2.03 & 2.03 & 2.00 & 2.00 & 2.00 & & & 2.00 & 2.04 & 2.07 & 2.05 & 2.06 \\
       & Sd & 0.60 & 0.40 & 0.27 & 0.17 & 0.03 & & & 0.52 & 0.41 & 0.26 & 0.16 & 0.04 \\
       & $\Delta Z_{(m)}$ & 0.23 & 0.23 & 0.23 & 0.22 & 0.23 & & & 0.06 & 0.08 & 0.07 & 0.07 & 0.07 \\
       & Sd, $\Delta Z_{(m)}$ & 0.03 & 0.02 & 0.01 & 0.01 & 0.01 & & & 0.01 & 0.01 & 0.01 & 0.00 & 0.00 \\
       & time & 0.17 & 0.27 & 0.63 & 1.28 & 16.25 & & & 1.90 & 2.42 & 3.22 & 4.32 & 22.94 \\
    \bottomrule
    \bottomrule
  \end{tabular}
  \label{logCauchyTable}
\end{table}

% Table generated by Excel2LaTeX from sheet 'GMMToy'
\begin{table}[hbtp]
  \centering
  \scriptsize
  \caption{\scriptsize Results for Simulation 3 (GMM, number of group $K$ = 3). $M = 1$ means single observation per individual and $N$ is number of individuals. $\pi_k$ is Multi-nomial distribution parameters of mixture, with Gaussian distribution parameter $\mu_k$ and $\sigma^2_k$. ``Est" and ``Sd" represent mean and standard deviation of estimator $\widehat{\btheta}$ over 5000 Monte Carlo experiments. Estimations are evaluated by average error $\widehat{\btheta} - \btheta$ in ``bias", and the ratio of correct prediction in ``Accuracy".}
    \begin{tabular}{ccrrrrrrrrrr}
    \toprule
    \toprule
    \multicolumn{12}{c}{$N$ = 50} \\
    
    \multicolumn{2}{c}{\multirow{2}[1]{*}{True Value}} & \multicolumn{1}{l}{$\mu_1$} & \multicolumn{1}{l}{$\mu_2$} & \multicolumn{1}{l}{$\mu_3$} & \multicolumn{1}{l}{$\sigma_1^2$} & \multicolumn{1}{l}{$\sigma_2^2$} & \multicolumn{1}{l}{$\sigma_3^2$} & \multicolumn{1}{l}{$\pi_1$} & \multicolumn{1}{l}{$\pi_2$} & \multicolumn{1}{l}{$\pi_3$} & \multicolumn{1}{l}{Accuracy} \\
    \multicolumn{2}{c}{} & -3.00 & 0.00 & 3.00 & 1.00 & 1.00 & 1.00 & 0.30 & 0.50 & 0.20 & \% \\
    
    \multirow{3}[0]{*}{EM} & bias  & 0.20 & -0.01 & -0.68 & -0.05 & 0.26 & 0.03 & -0.01 & -0.02 & 0.02 &  \\
          & Est   & -2.81 & -0.01 & 2.32 & 0.95 & 1.26 & 1.03 & 0.29 & 0.49 & 0.22 & 80.51 \\
          & Sd    & 1.22 & 0.58 & 1.93 & 0.88 & 1.10 & 1.11 & 0.14 & 0.18 & 0.13 & 11.20 \\
    \multirow{3}[1]{*}{MILE} & bias  & 0.08 & 0.03 & -0.53 & -0.34 & 0.03 & -0.39 & -0.01 & 0.01 & 0.00 &  \\
          & Est   & -2.92 & 0.03 & 2.47 & 0.66 & 1.03 & 0.61 & 0.29 & 0.51 & 0.20 & 80.90 \\
          & Sd    & 1.23 & 0.62 & 2.03 & 0.59 & 0.97 & 0.68 & 0.13 & 0.17 & 0.11 & 11.23 \\
    
    \multicolumn{12}{c}{$N$ = 500} \\
    
    \multirow{3}[0]{*}{EM} & bias  & -0.03 & -0.03 & -0.24 & 0.03 & 0.22 & 0.23 & -0.01 & -0.02 & 0.02 &  \\
          & Est   & -3.03 & -0.03 & 2.77 & 1.03 & 1.22 & 1.23 & 0.29 & 0.49 & 0.22 & 86.98\\
          & Sd    & 0.28 & 0.20 & 0.83 & 0.36 & 0.92 & 0.85 & 0.07 & 0.08 & 0.07 & 6.03\\
    \multirow{3}[1]{*}{MILE} & bias  & -0.12 & -0.03 & 0.12 & -0.21 & -0.14 & -0.29 & -0.01 & 0.01 & -0.00 &  \\
          & Est   & -3.12 & -0.03 & 3.12 & 0.79 & 0.86 & 0.71 & 0.29 & 0.51 & 0.20 & 87.99\\
          & Sd    & 0.24 & 0.17 & 0.26 & 0.24 & 0.63 & 0.21 & 0.07 & 0.09 & 0.04 & 4.85\\
    \bottomrule
    \bottomrule
    \end{tabular}%
  \label{GMMTable}%
\end{table}%

\begin{table}[hbtp]
  \centering
  \scriptsize
  \caption{\scriptsize Results of Simulation 4 (Bayesian Segmented Regression). $M$ is numbers of observations per individual and $N$ is number of individuals. ``Est" and ``Sd" represent mean and standard deviation of estimator $\widehat{\btheta}$ over 5000 Monte Carlo experiments. Latent variable estimation is evaluated by the average bias $\Delta Z = mean(\widehat{Z} - Z)$.}
  \begin{tabular}{clrrrrrrclrrrrrr}
    \toprule
    \toprule
    \multicolumn{8}{c}{$N$ = 4} & \multicolumn{8}{c}{$N$ = 25} \\
    \multirow{2}[2]{*}{$M$} & \multicolumn{1}{c}{\multirow{2}[2]{*}{True}} & \multicolumn{1}{l}{$\alpha$} & \multicolumn{1}{l}{$\beta$} & \multicolumn{1}{l}{$\beta_1$} & \multicolumn{1}{l}{$\beta_2$} & \multicolumn{1}{l}{$a$} & \multicolumn{1}{l}{$\Delta Z$} 
    & \multirow{2}[2]{*}{$M$} & \multicolumn{1}{c}{\multirow{2}[2]{*}{True}} & \multicolumn{1}{l}{$\alpha$} & \multicolumn{1}{l}{$\beta$} & \multicolumn{1}{l}{$\beta_1$} & \multicolumn{1}{l}{$\beta_2$} & \multicolumn{1}{l}{$a$} & \multicolumn{1}{l}{$\Delta Z$} \\
      & & 5.00 & 5.00 & 1.00 & -1.00 & 1.00 & 0.00 
        & & & 5.00 & 5.00 & 1.00 & -1.00 & 1.00 & 0.00 \\
    \multirow{2}[1]{*}{20} & Est & 21.61 & 23.65 & 1.81 & -1.70 & 1.10 & 0.07 
    & \multirow{2}[1]{*}{20} & Est & 10.43 & 10.63 & 0.90 & -0.91 & 0.86 & -0.01 \\
    & Sd  & 26.28 & 29.71 & 3.83 & 4.01  & 0.29 & 0.78 
    &     & Sd  & 3.31  & 4.03  & 0.19 & 0.20  & 0.12 & 0.89 \\
    \multirow{2}[0]{*}{40} & Est & 21.53 & 22.26 & 1.09 & -1.10 & 1.07 & -0.01 
    & \multirow{2}[0]{*}{40} & Est & 9.29  & 9.20  & 0.87 & -0.91 & 0.86 & -0.06 \\
    & Sd  & 23.48 & 24.20 & 0.27 & 0.30  & 0.17 & 0.37 
    &     & Sd  & 3.00  & 3.47  & 0.16 & 0.17  & 0.09 & 0.78 \\
    \multirow{2}[0]{*}{60} & Est & 17.60 & 17.33 & 1.05 & -1.06 & 1.03 & 0.01 
    & \multirow{2}[0]{*}{60} & Est & 8.70  & 8.64  & 0.86 & -0.91 & 0.85 & -0.05 \\
    & Sd  & 20.16 & 20.27 & 0.19 & 0.21  & 0.15 & 0.27 
    &     & Sd  & 2.54  & 3.12  & 0.15 & 0.16  & 0.07 & 0.74 \\
    \multirow{2}[1]{*}{200} & Est & 15.02 & 14.86 & 1.02 & -1.00 & 1.00 & 0.02 
    & \multirow{2}[1]{*}{200} & Est & 8.08  & 7.94  & 0.86 & -0.91 & 0.87 & -0.06 \\
    & Sd  & 18.34 & 17.40 & 0.09 & 0.09  & 0.07 & 0.16 
    &     & Sd  & 2.25  & 2.57  & 0.12 & 0.12  & 0.04 & 0.63 \\
    \bottomrule
    \bottomrule
  \end{tabular}
  \label{BSRTable}
\end{table}

%\section{Acknowledgments}
%We thank the anonymous reviewers for their valuable suggestions and efforts.

\bibliographystyle{chicago}
\bibliography{reference}

\end{document}